\documentclass[letterpaper, paper, 10pt]{AAS}		

\usepackage{bm}
\usepackage{amsmath}
\usepackage{siunitx} 
\usepackage[colorlinks=true, pdfstartview=FitV, linkcolor=black, citecolor= black, urlcolor= black]{hyperref}
\usepackage{overcite}
\usepackage{footnpag}			      	
\usepackage{gensymb}
\DeclareMathAlphabet{\mathpzc}{OT1}{pzc}{m}{it}
\usepackage{mathtools}
\usepackage{adjustbox}
\usepackage{graphicx}
\usepackage{array}
\usepackage{caption}
\usepackage{subcaption}
\usepackage{hyperref}
\usepackage{amssymb}
\usepackage[super]{nth}
\usepackage{color}
\usepackage{soul}
\usepackage[bottom]{footmisc}
\usepackage{tabularx}
\usepackage[table,xcdraw]{xcolor}
\usepackage{ragged2e} 
\DeclareMathOperator*{\argmin}{arg\,min}
\DeclareMathOperator*{\argmax}{arg\,max}

\usepackage{siunitx}

\def\bfg{\mbox{\boldmath $g$}}

\def\bfr{\mbox{\boldmath $r$}}

\def\bfv{\mbox{\boldmath $v$}}

\def\bflambda{\mbox{\boldmath $\lambda$}}

\PaperNumber{19-894}

\title{GTOC X: SOLUTION APPROACH OF TEAM SAPIENZA-POLITO}
\author{Alessandro Zavoli\footnote{Research Assistant, Department of Mechanical and Aerospace Engineering, Sapienza University of Rome, Via Eudossiana 18 - 00184, Rome, Italy. alessandro.zavoli@uniroma1.it}
, Lorenzo Federici\footnote{PhD student, Department of Mechanical and Aerospace Engineering, Sapienza University of Rome, Via Eudossiana 18 - 00184, Rome, Italy. lorenzo.federici@uniroma1.it}, Boris Benedikter\footnote{PhD student, Department of Mechanical and Aerospace Engineering, Sapienza University of Rome, Via Eudossiana 18 - 00184, Rome, Italy. boris.benedikter@uniroma1.it}, 
\\Lorenzo Casalino\footnote{Full Professor, Department of Mechanical and Aerospace Engineering, Turin Polytechnic, Corso Duca degli Abruzzi 24 - 10129, Turin, Italy. lorenzo.casalino@polito.it}, Guido Colasurdo\footnote{Full Professor, Department of Mechanical and Aerospace Engineering, Sapienza University of Rome, Via Eudossiana 18 - 00184, Rome, Italy. guido.colasurdo@uniroma1.it}
}

\begin{document}

\maketitle

\begin{abstract}
     This paper summarizes the solution approach and the numerical methods developed by the joint team Sapienza University of Rome and Politecnico di Torino (Team Sapienza-PoliTo) in the context of the  \nth{10} Global Trajectory Optimization Competition.
  The proposed method is based on a preliminary partition of the galaxy into several small zones of interest, where partial settlement trees are developed, in order to match   a (theoretical) optimal star distribution.
  A multi-settler stochastic Beam Best-First Search, 
  that exploits a guided  multi-star multi-vessel transition logic, is proposed for solving a coverage problem, where the number of stars to capture and their distribution within a zone is assigned.
  The star-to-star transfers were then optimized through an indirect procedure.
  A number of refinements, involving \textit{settle time re-optimization}, \textit{explosion}, and \textit{pruning}, were also investigated.
  The submitted 1013-star solution, as well as an enhanced 1200-point rework, are presented.

\end{abstract}

\section{Introduction}
This paper provides an overview of the solution approach and the numerical methods developed by the joint team Sapienza University of Rome and Politecnico di Torino (Team Sapienza-PoliTo) in the context of the \nth{10} Global Trajectory Optimization Competition, also known as GTOC X.

The proposed problem, named ``Settlers of the Galaxy'', 
represents a unique challenge, 
where participants are asked to design
a settlement tree for colonizing the galaxy, \emph{i.e.}, a pool of one hundred thousand candidate stars,
scoring points according to a prescribed merit function $J=J_1 J_2 J_3$, that rewards 
large settlements, uniformly distributed in space ($J_2$) and  an efficient use of propulsion ($J_3$).
Solutions submitted earlier in the competition are also slightly rewarded ($J_1$).
More precisely:
\begin{align}
J_2 &= \frac{N}{1+10^{-4}\cdot N\left(E_r+E_\theta\right)} \\
J_3 &= \frac{\Delta V_{max}}{\Delta V_{used}} 
\end{align}
being $N$ the total number of settled stars, $E_r$ and $E_\theta$ the uniformity errors in radius $r$ and final polar angle $\theta_f$, respectively, $\Delta V_{used}$ the total $\Delta V$ employed and $\Delta V_{max}$ the maximum admissible $\Delta V$.

A set of rules is provided, which constrains the number of starting vessels, namely Mother Ships and Fast Ships, as well as the number of new vessels, namely Settler Ships, that may depart from any settled star.
Maneuverability limits for each type of vessel, \emph{i.e.}, maximum number and magnitude of allowed velocity increments, are also assigned.
A peculiar, non-keplerian, central-force gravitational model 
is assumed for the motion of both stars and ships, 
inspired by the circular motion observed for actual stars in the Milky Way.
Further details about the GTOCX problem are reported in the statement paper.\cite{GTOCXProblem}

The paper  is organized as follows.
First,
an overview of the global strategy pursued by this team is presented.
An indirect method for the ships trajectory optimization is summarized, together with an approximate solution based on the impulsive linear rendezvous problem.
A set covering problem is defined, with the aim of attaining an \textit{a priori} chosen star distribution over some zones of the galaxy.
A  multi-settler stochastic Beam Best-First Search, that exploits a guided multi-star multi-vessel transition logic, is used for solving this novel problem.  A number of refinements are then discussed. The submitted 1013-star solution, as well as an enhanced version with a score above 1200 points, are presented.

\section{Strategy Overview}

According to the experience acquired in the past GTOC editions,\cite{GTOC5Casalino,GTOC6,GTOC8}
the solution approach of Team Sapienza-PoliTo mostly depended on the identified problem features.
On the basis of a preliminary analysis of the merit function $J$, 
the minimization of the error functions $E_r$ and $E_\theta$ was deemed essential for maximizing the score, hence it was assumed as main design goal.
Indeed, one should notice that if the error functions $E_r$ and $E_\theta$ are small, $J$ increases almost linearly with the number of settled stars $N$.

A problem reformulation was devised in order not to deal directly with the error functions whilst providing an almost equivalent, yet more practical, way to handle the issue.
The galaxy  was partitioned according to a $30\times 32$ grid, defined over the ($r$, $\theta_f$) plane, where the $k_r$-th horizontal row (or ``ring'') was defined by $r\in \left[1 + k_r,\, 2 + k_r \right]$,  with $k_r \in [1,\,30]$, and the $k_\theta$-th column (or   ``slice'') by  $\theta_f \in \left[-\pi + \frac{\pi}{16} (k_\theta-1),\, -\pi + \frac{\pi}{16} k_\theta \right]$, with $k_\theta \in [1,\,32]$.

Next, an optimal distribution of stars in each row and column for minimizing the error functions $E_r$ and $E_\theta$  was searched for,
by assuming the center of each grid cell as a proxy for the actual location of each star inside that cell.
This assumption was motivated by the fact that if the star distribution is sufficiently dense and the number of captured stars is sufficiently large, 
then the local effects would be minimal.

As a result of this analysis, it was found out that the errors could be kept to a minimum by a settlement tree of $512 \times n$ stars, equally-distributed over the 32 columns and linearly-distributed over the 30 rows, with more settled stars at the higher radii than at lower ones.
Given the limited team experience in dealing with very large data structures, as well as the short available time,  a decision was made to limit the number of settled stars to about 1020, postponing all the possible solution refinements to the last few days of the competition.
Randomly sampling stars according to the proposed distribution provided a (reasonable-enough) empirical proof for justifying the proposed problem reformulation.
Numerical investigations also suggested that small deviations from the target distribution leading to errors smaller than unity could be accepted without severely hindering the overall score.
Conversely, an unsupervised increment of the number of stars would have resulted in a significant score drop.

Given the small number of vessels at Year Zero,
the role of Settler Ship multiplication in creating the settlement tree was apparent.
Mobility charts were drawn for the Settler Ships in order to provide a graphical tool for assessing the maneuvering capability of a vessel, that is, the location of all reachable arrival stars in the ($r$, $\theta_f$) plane as a function of departure radius, flight time, and departure epoch.

A spiral occupation was recognized as the only viable option, being a global uniform coverage of the space unfeasible due to the strict time constraints. 
Settler Ship mobility charts were used for guiding the search of the opening maneuvers, that is,  Mother Ships  and Fast Ships trajectories departing from Sol at Year Zero.
In order to keep the search for the opening trajectories  simple enough, 
Mother Ship trajectories intercepting more than three stars had not been investigated.
The target number (1020) of settled stars  was so low that 
two or three ``root'' settlements per Mother Ship were deemed sufficient for attaining the desired number of settled stars. 
Opening maneuvers had been concurrently optimized in order to provide ``root'' settlements in galaxy locations sufficiently far apart  to allow 
a proper ``growth'' of the spiral settlement tree,
without having multiple ``sub-trees'' originating from different roots competing each other for the same stars.
Eventually, 9 non-overlapping zones, that is, sets of contiguous cells of the ($r$, $\theta_f$) grid, had been identified.
In each zone a set-covering problem is posed
as
finding the optimal settlement tree, for an assigned
departure (or ``root'') star and epoch,
that allows to settle a prescribed number of stars per cell within the zone.

A multi-settler stochastic Beam Best-First Search (BBFS) has been developed for solving this unusual problem.
The best-first variant was motivated by the need to attain an answer to each set-covering problem as soon as possible, making possible to iteratively adjust all the partial settlement trees. 
Settlement trees of up to 100 stars were generated in a few minutes on a standard laptop.

A settlement tree with $1013$ stars was found at the end of this first phase. 
A Genetic Algorithm was adopted to slightly refine the transfer times, leading to the submitted solution. 
An enhanced solution with $1220$ stars and 
an overall score above 1200 points was also investigated by adopting a greedy strategy aimed at increasing the number of settled stars (hence the cost factor $J_3$)
without destroying the spatial distribution (hence with almost constant $J_2$) and then removing unnecessary colonies.
Unfortunately, the competition ended before we were able to submit it.

\section{Trajectory optimization in the Milky Way}

\subsection{Dynamical Model}
A non-Keplerian dynamical model governs the motion of both stars and ships.
A central-force law is assumed, which approximately models the circular motion observed for actual stars in our Milky Way.
The equations of motion are formulated as:
\begin{align}
    \dot{\bfr} &= \bfv \label{eq:nonlin-r}\\
    \dot{\bfv} &= \bfg = - \frac{\bfr}{r} f(r) \label{eq:nonlin-v}
\end{align}
where $\bfr$ and $\bfv$ are, respectively, the position and the velocity vector, 
and $\bfg$ is the acceleration vector directed towards the galactic center, of modulus $f(r) = v_c^2/r$, with 
$v_c$ the circular orbit speed for a body at radius $r$ from the galactic center, expressed as the inverse of an $\nth{8}$-order polynomial function of known coefficients $k_0, \ldots, k_8$, that is:
\begin{equation}
\displaystyle
    v_c(r) = \frac{1}{\sum\limits_{i=0}^{8}{k_i r^i}}
    \label{eq:circvel}
\end{equation}


\subsection{Indirect Optimization}

Only impulsive maneuvers are admitted. According to the vessel type, up to five impulses can be performed.
The ``basic'' maneuver is a two-impulse transfer between two stars with fixed departure and arrival time. The problem at hand is equivalent to Lambert's problem in Keplerian dynamics: the initial velocity components (three unknowns) must be determined in order to match the target position at the arrival time (three boundary conditions). The problem has a single solution and is not subject to optimization. It is worth noting that the trajectories improve when the time of flight is increased, as the competition overall time constraint forces times of flight that are much shorter than the optimal one (e.g., the time competing to the equivalent of a Hohmann transfer).

Despite their simplicity, two-impulse maneuvers are usually optimal or close to optimality. However, impulses may be larger than allowable limits, either on the single impulse magnitude and/or on the overall $\Delta V$ value. 
Trajectories exceeding limits must be discarded for Fast Ships, which can perform at most two impulses.
Instead, more impulses are available for Mother Ships (up to 3) and Settler Ships (up to 5), and multiple-impulse trajectories are sought when necessary.

A multiple-impulse trajectory requires the definition of an optimal control problem that can be solved through an indirect approach\cite{Bryson1979,ZavoliAlaska,Casalino1999}.
For the dynamical system of Eqs.~\eqref{eq:nonlin-r}-\eqref{eq:nonlin-v}, 
the Hamiltonian $H = \bflambda_{\bfr}^T \bfv + \bflambda_{\bfv}^T \bfg$ is defined and the Euler-Lagrange equations are derived for the adjoint variables:
\begin{equation}
\dot \bflambda_{\bfr}^T = -\frac{\partial H}{\partial \bfr}  \hspace{2cm} \dot \bflambda_{\bfv}^T = -\frac{\partial H}{\partial \bfv}
\end{equation}
The problem is completed by the boundary conditions on state variables, that must ensure to reach the target and fulfill the $\Delta V$ limitations, and the boundary conditions for optimality. 

For Settler Ships, up to 4 impulses are considered, since no trajectory requiring 5 impulses has been found. The sum of the impulses magnitude is minimized; constraints on the initial and/or final impulse magnitude are introduced when limits are exceeded. The boundary conditions for optimality require impulses parallel to the velocity adjoint vector (primer vector) and the same primer magnitude at non-constrained impulses.\cite{Casalino1998} The primer magnitude must be maximum at intermediate impulses when the time is not constrained; however, this would lead to a violation of the minimum timing between maneuvers and the corresponding times must therefore be constrained.  

Although Mother Ships can have encounters with multiple stars, they are treated with a similar approach to minimize the overall $\Delta V$. Two classes of trajectories have been exploited in the solution: impulse-intercept-impulse-intercept-impulse-intercept and impulse-impulse-intercept-impulse-intercept. Again constraints on $\Delta V$s are introduced when limits are exceeded. Impulses are always at the minimum allowable time. Intercept times are left free and the conditions for optimality requires that the scalar product of primer and Settlement Pod $\Delta V$ is zero.




\subsection{Approximate Solution: the Minimum-Time Linear Impulsive Rendezvous Problem}

A proper optimization of each vessel trajectory  
is mandatory, as saving $\Delta V$ positively impact on $J_3$, hence on the total score.
However, the preliminary definition of the (partial) settlement trees
requires the evaluation and the comparison of a possibly enormous amount of transfers connecting any pair of stars, with different departure and arrival epochs.
The need for a fast computational tool, able to provide 
a rough evaluation of transfer time and cost, is thus apparent.

The competition overall time constraint dominates the problem;
in order to achieve a final ``populous'' settlement tree,
short (and fast) trajectories between relatively close stars are very likely.
This considerations suggest the use of 
the solution of an impulsive linearized rendezvous problem
as a way for rapidly estimate the transfer $\Delta V$s as a function of transfer duration.

Linearized relative motion, 
that is, Hill-Clohessy-Wiltshire equations\cite{wh1960terminal}, 
are attained by assuming as reference body the target star, which moves on a circular orbit of radius $\tilde r_f$ with angular velocity $\tilde \omega_f = v_c(\tilde r_f)/\tilde r_f$.
The natural (linearized) motion of a vessel, relative to the target,
can be written in compact matrix form as:
\begin{align}
    \begin{bmatrix} \delta \bfr \\ \delta \bfv \end{bmatrix} = 
    \begin{bmatrix}
    M(\tau) & N(\tau) \\
    S(\tau) & T(\tau)
    \end{bmatrix}
    \begin{bmatrix} \delta \bfr_0 \\ \delta \bfv_0 \end{bmatrix}
\label{eq:hill-STM}
\end{align}
where $\delta \bfr$, $\delta \bfv $ are the position and velocity vector of the vessel in the relative reference frame, respectively, $\tau = (t_f-t_0)$ is the transfer time, and $M(\tau)$, $N(\tau)$, $S(\tau)$, $T(\tau)$ are $3\times 3$ matrices composing the state transition matrix, here reported for the sake of completeness:
\small
\begin{align}
M &= \begin{bmatrix} 
    4-3\cos{\tilde \omega_f\tau} &0 &0  \\
 6(\sin{\tilde \omega_f\tau} - \tilde \omega_f\tau) & 1& 0\\
 0&0 &\cos{\tilde \omega_f\tau} \end{bmatrix} &
N &= \begin{bmatrix} 
\frac{1}{\tilde \omega_f} \sin{\tilde \omega_f\tau}&  \frac{2}{\tilde \omega_f}(1-\cos{\tilde \omega_f\tau}) & 0 \\
 -\frac{2}{\tilde \omega_f}(1-\cos{\tilde \omega_f\tau})& \frac{4}{\tilde \omega_f}\sin{\tilde \omega_f\tau}-3\tau& 0\\
 0& 0 & \frac{1}{\tilde \omega_f} \sin{\tilde \omega_f\tau}\end{bmatrix}\\
S &= \begin{bmatrix} 
3\tilde \omega_f\sin{\tilde \omega_f\tau} & 0 &0 \\
-6\tilde \omega_f\sin{\tilde \omega_f\tau}  &0 & 0\\
0 &0 & -\tilde \omega_f\sin{\tilde \omega_f\tau}\end{bmatrix}
 &
T &= \begin{bmatrix} 
\cos{\tilde \omega_f\tau} & 2\sin{\tilde \omega_f\tau} & 0\\
-2\sin{\tilde \omega_f\tau} &4\cos{\tilde \omega_f\tau}-3 &0 \\
0 & 0 & \cos{\tilde \omega_f\tau} \end{bmatrix}
\end{align}
\normalsize

Optimal impulsive linear rendezvous problems have been deeply studied\cite{prussing1967optimal}.
Two-impulse solutions are particularly simple to compute, as no iterative procedure is involved; only a few matrix operations are required. 
%
The rendezvous condition requires null relative displacement and velocity at the final time, that is $\delta \bfr_f$ = 0 and  $\delta \bfv_f = 0$.
Let $\delta \bfr_0$, $\delta \bfv_0$ be known initial conditions. 
The departure velocity $\delta \bfv_0^+$ required to intercept the target at time $\tau$ and the corresponding velocity at the intercept $\delta\bfv_f^{-}$ are 
evaluated as:
\begin{align}
    \delta \bfv_0^+ &= -N^{-1}(\tau)\left( M(\tau) \delta \bfr_0\right) \label{eq:dv0+}\\
    \delta\bfv_f^{-} &= S(\tau) \delta \bfr_0 + T(\tau) \delta \bfv_0^+ \label{eq:dvf-}
\end{align}
The total transfer cost is readily computed as :
%
\begin{equation}
    \Delta V_{tot} = \Delta V_0 + \Delta V_f = \left\| \delta\bfv_0^{+}  - \delta\bfv_0 \right\|  + \left\| \delta\bfv_f^{-}\right\|     \label{eq:DVhill}
\end{equation}

%
%

Discrepancies between full-dynamics and linear solutions were minimal in most cases of interest (below $2\div 5$ km/s).
Indeed, as for the full-dynamics problem, impulse magnitude limits are likely to be violated.
Therefore, time-minimum two-impulse transfer that satisfy the constraints on both single and overall impulse magnitude have been searched for.
A simple line search on a time grid with 1 Myr step provides reasonable values of the minimum admissible rendezvous time $\tau$ for any assigned pair of stars and departure epoch.

Eventually, given a starting star $s$ and a starting epoch $t_s$,
a minimum-time Hill neighborhood $\mathcal{N}_s = \mathcal{N}(s,t_s)$ has been defined as the set of all stars that can be reached in the available time window departing from $s$ at time $t_s$, while respecting the vessel $\Delta V$ limits
on both single impulse magnitude ($\Delta V_i < \Delta V_i^{max}$) and their sum ($\Delta V_{tot} < \Delta V_{tot}^{max}$).
In order to account for possible improvements that would inevitably come with the  full-optimization procedure,
$\Delta V$ constraints are (moderately) relaxed, with the tolerance increasing together with the flight time.
Minimum-time hill neighborhoods have been extensively used in the construction of partial settlements trees.
The adopted settings, based on a trial and error procedure, are summarized in Table~\ref{tab:DVtol}.

\begin{table}[!h]
\centering
\caption{Acceptance thresholds for the generation of minimum-time Hill neighborhoods.  }
\label{tab:DVtol}
\begin{tabular}{lcc}
    \hline\hline
   &   \rule{0pt}{12pt} $\Delta V_i^{max}$ [km/s] & $\Delta V_{tot}^{max}$ [km/s] \\[2pt]    
       \hline
    $\tau<~4$ Myr & 170 & 340  \\
    \phantom{$\tau_f$} $<~9$ Myr & 175& 350\\
    \phantom{$\tau_f$} $< 15$ Myr & 190& 360 \\
    \multicolumn{1}{c}{otherwise}      &300 & 400 \\
    \hline\hline
\end{tabular}
\end{table}



\section{Searching for Settlement Trees} \label{sez:search}



A major role in the problem solution is played by the {creation of}  (partial) settlements trees,
that is, a web of vessel trajectories connecting pairs of stars, with origin at one given  \textit{root} star, either settled by a Mother Ship or a Fast Ship.
The search for settlements trees is aimed at providing a proper spatial distribution,  that keeps the error functions sufficiently low, 
while ensuring the feasibility of all transfers.
To attain this result,
a number of \textit{zones}, that is, sets of contiguous cells
of the $30\times 32$ grid  defined  over  the  ($r$, $\theta_f$)  plane,
are identified.
Each zone $Z$ encompasses a {prescribed} number of adjacent rows $K_r^{Z} = \{ k_{r,min}^Z, \dots, k_{r,max}^Z\}$ and columns $K_\theta^{Z} = \{k_{\theta,min}^Z, \dots, k_{\theta,max}^Z\}$.


In each zone $Z$, a \textit{set-covering problem} is issued.
%
%
The goal is to find, if it exists, a
settlement tree,
starting at a given root star and epoch,
that realizes a spatial occupancy as close as possible to 
a target one.
%
%
The target \textit{occupancy distribution} in a zone could be provided either on a cell-by-cell basis (i.e., by imposing the number of stars to settle in each cell), or on a row and column basis (if the stars to settle are imposed in each row and column).
%
%
Having in mind the (preliminary studied) global $1020$-star occupancy, 
the row/column target formulation appears simpler to devise and implement.
Unfortunately, this formulation results in a very broad, hence extremely slow, search. Instead, since a same row/column distribution can be achieved with several (different) cell-by-cell distributions, 
the use of a formulation by cells narrows down the search, speeding up the solution process.
In addition, a cell-by-cell distribution allows to completely remove all transfers ending to stars belonging to a cell where the occupancy requirement has already be filled (both in case the target occupancy is zero, or the target number of stars in that cell has already been settled), so simplifying further the search. 
%
%
These reasons led us to prefer an objective distribution by cells instead of by rows and columns, in each of the identified zones.
As a minor drawback, 
a few iterations over the target distribution matrix of each zone were needed,
in order to produce an actual overall distribution complying with the $1020$-star desired row/column occupancy.

%


 \subsection{Set covering problem definition}

The set-covering problem poses as a ``standard'' \textit{search problem} \cite{russell2016artificial}, where the objective is to look for a sequence of \textit{actions}, \emph{i.e.}, transitions from a state to the next one (or \textit{successor}), that, starting from a known initial state, reaches a \textit{goal} state, that is, an optimal state according to a given \textit{performance measure}.

The set-covering problem can be formally defined as a search problem, as soon as the following components are specified: (i) the state space of the problem, that is, the set of all states reachable from the initial state by any sequence of actions, (ii) a performance measure, to compare the states each other and push the search towards the goal, and (iii) the transition model, which allows to generate the successors of a given state.

\subsubsection{State Space.} \label{sez:states}
Each state in the state space represents a settlement tree, \emph{i.e.}, an un-ordered list of stars $\mathcal{S} =\{ s_0, s_1, \dots, s_n\}$ settled at times $t_0, t_1, \dots t_n$. The k-th star $s_k$ is labeled with information about: i) the star ID, $ID_k$, ii) the ``parent'' star, $s_p$, where the in-bound Settler Ship comes from (with $p \in [0,n-1]$), iii) a maximum of three ``offspring'' stars, $\mathcal{O}_k = \bigcup_{i = 1}^{l \leq 3} s_{o_i}$, reached by the out-bound Settler Ships (with $o_i \in [1,n]\,\, \forall i$), and iv) the star settlement epoch $t_k$.
Since only the minimum (2-Myr) waiting time has been considered at each settled star, the evaluation of the transfer time from a parent star to any of its offspring is straightforward, once the settlement epochs are known.
The initial state of the search procedure is the settlement tree made up only by the root star $s_0$, settled at a known epoch $t_0$.
\subsubsection{Performance Measure.} \label{sez:cost}
The explored states are ranked according to the following cost function (to be minimized):
\begin{equation}
    \phi = \|H - G\| + \frac{t_n}{100}
    \label{eq:g}
\end{equation}
%
$H$ and $G$ are two matrices such that, for each cell ($k_r,\,k_\theta$) composing the zone $Z$ under search, $H(k_r,\,k_\theta)$ (respectively, $G(k_r,\,k_\theta)$) reports the actual (respectively, target) number of stars settled (respectively, to settle) in that cell, and $t_n$ is the final settling time in Myr.
Being $\|H - G\|$ an integer value and $\frac{t_n}{100} < 1$, the final settling time becomes important only once the target distribution has been achieved.

\begin{figure} [h!]
\captionsetup[subfigure]{justification=centering}
    \centering
    \begin{subfigure}{0.45\textwidth}
    \centering
       \includegraphics[width=\textwidth]{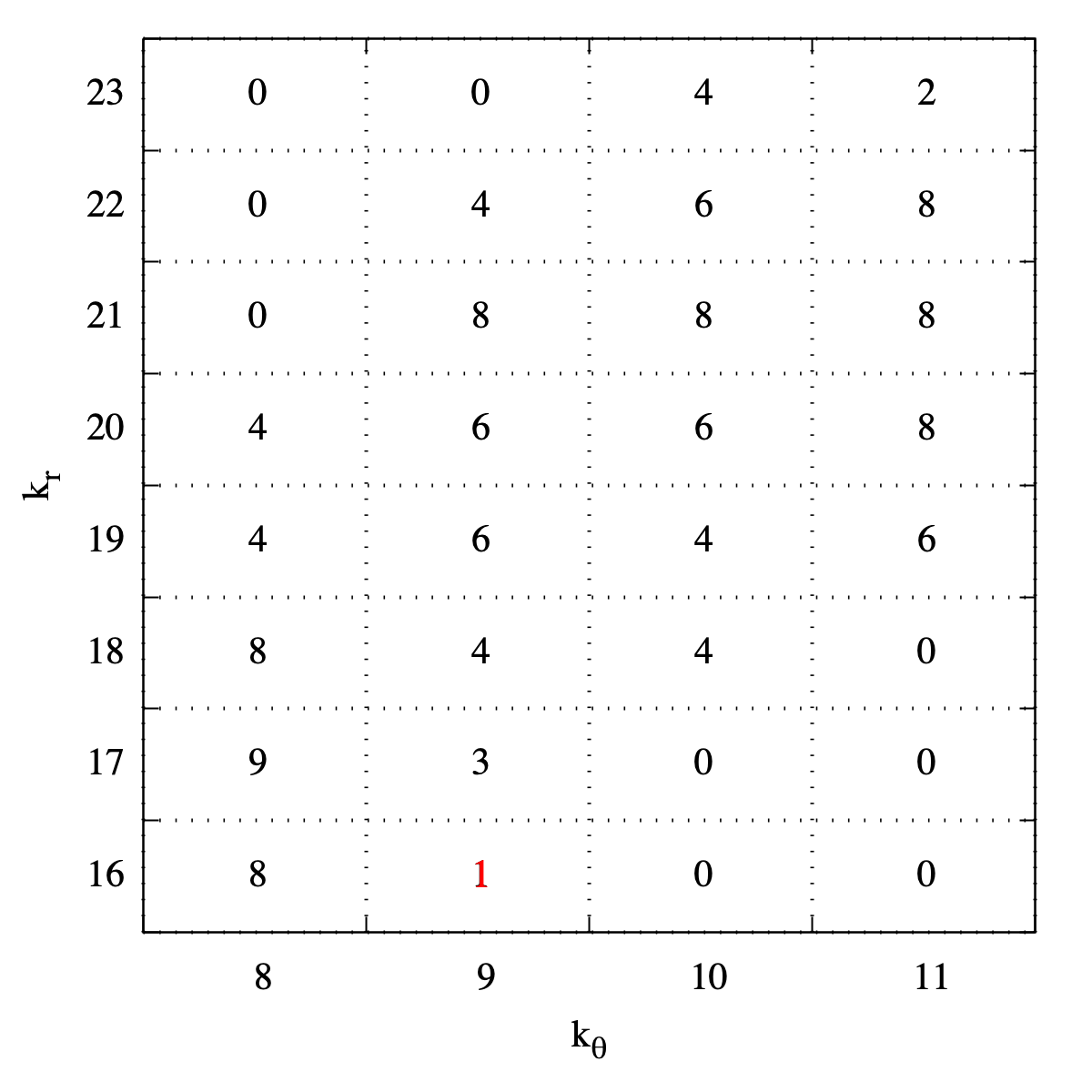}
       \caption{Distribution matrix $G$. The red number indicates the root star.}
    \label{fig:M11mat}
    \end{subfigure}
    \begin{subfigure}{0.45\textwidth}
    \centering
       \includegraphics[width=\textwidth]{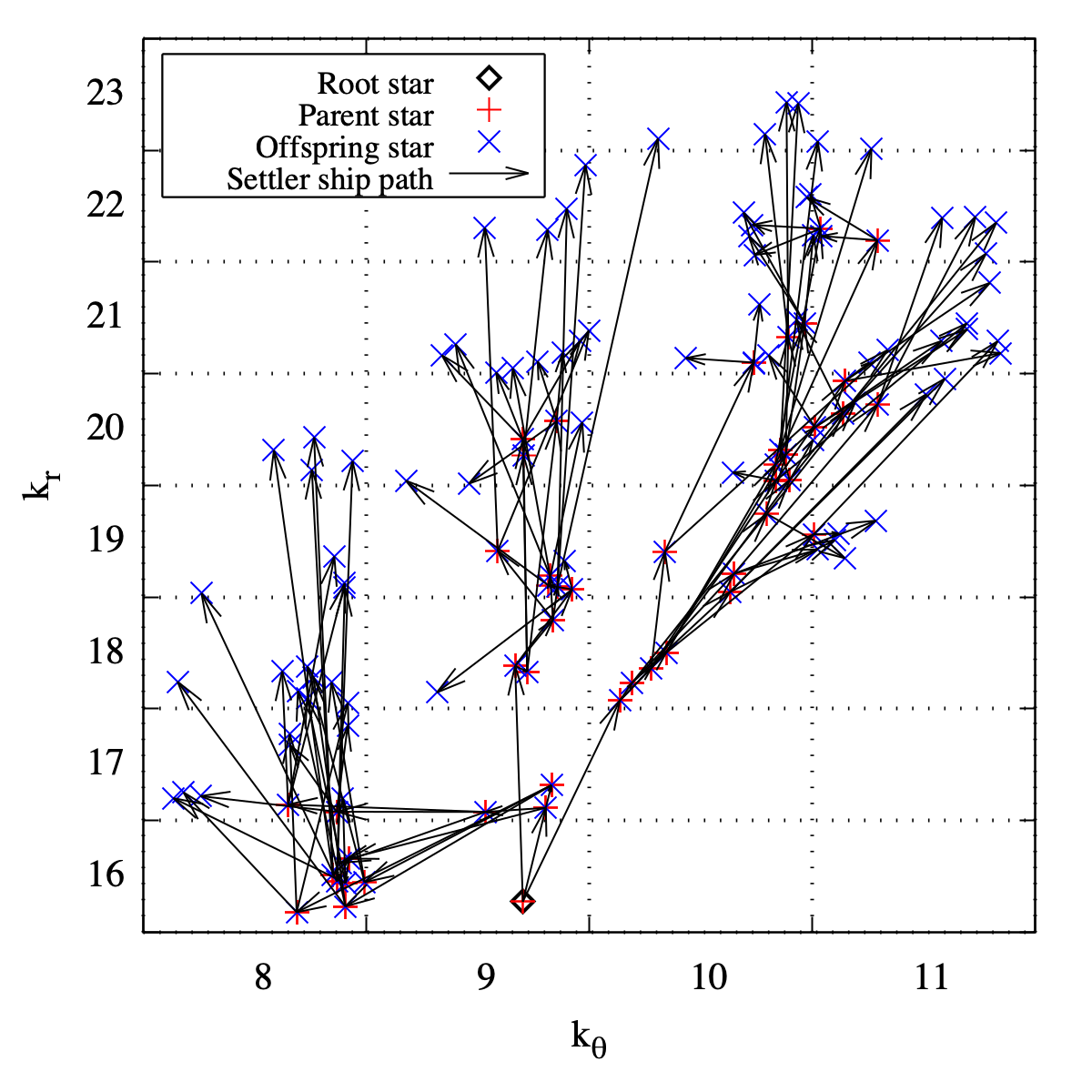}
        \caption{Settlement tree obtained as a result of the search procedure.}
    \label{fig:M11tree}
    \end{subfigure}
   \caption{Zone $M11$: target distribution (a), correspondent settlement tree (b).}
\end{figure}

A visual representation in the ($r$, $\theta_f$) plane of the best settlement tree (\emph{i.e.}, state) obtained in the zone labeled $M11$, 
with target distribution matrix $G$ reported in Figure~\ref{fig:M11mat},
is proposed in Figure~\ref{fig:M11tree}, where time information is neglected for the sake of clarity.

\subsubsection{Transition Model.} \label{sez:successors}
New states, or \textit{successors}, are built upon previously evaluated ones by using an incremental mechanism, 
that is, by settling additional stars using the available vessels from any star composing the (already constructed) tree.

Let $s_k$ be a ``fertile'' settled star belonging to the \textit{ancestor} state $\mathcal{S}$, \emph{i.e.}, 
a settled star with at least one available vessel.
A reachable minimum-time neighborhood $\mathcal{N}^\star_k$ of $s_k$ is defined 
as the set of all (unsettled) stars, in the considered zone, reachable from $s_k$ in the least possible time which produces admissible $\Delta V$s.
%
In order to keep the computational time reasonably low,
the (true) neighborhood, obtained by solving the 
optimization problem involving, for each arrival star, the determination of a minimum-time, $\Delta V$-constrained, $n$-impulse transfer in the full dynamical model,
is replaced by the minimum-time Hill neighborhood $\mathcal{N}_k$,
where the two-impulse solution of the linear rendezvous problem is used.
Neighborhoods are computed offline, with departure and arrival epochs discretized over a uniform time-grid with a 1-Myr step size.
For the sake of simplicity, non-settled stars located in cells where 
current occupancy matches the target distribution $G$ are not included in any neighborhood, thus speeding up the computation.



Several transition models may be defined for the problem under investigation.
At a first attempt, a \textit{one-vessel} transition model was considered, that is,
any successor state $\mathcal{S}'$ is obtained by settling one additional star $s_o$ that belongs to the neighborhood of one of the $n_f$ fertile stars of the ancestor $\mathcal{S}$:

\begin{equation}
    \mathcal{S}' = \mathcal{S} \cup \left\{s_o \in \mathcal{N}_k \, |\, s_k \in \mathcal{S} \wedge |\mathcal{O}_k| <3 \right\}
\end{equation}

Up to $n_f\,|\mathcal{N}_k|$ different successor states per ancestor $\mathcal{S}$ may be constructed with this model,
being $|\mathcal{N}_k|$ the number of unsettled stars composing the minimum-time Hill neighborhood of a fertile star $s_k \in \mathcal{S}$. 



While particularly simple, this transition model was not really effective: indeed, the probability to obtain duplicate (successor) states starting from different ancestor states is quite high; as a consequence, the search algorithm has to keep in memory a great number of states very similar (if not equal) each other, and must perform a high number of steps (\emph{i.e.}, transitions) in the search space in order to reach the target distribution, thus resulting in a very slow procedure.

Therefore, 
a \textit{multi-star multi-vessel} transition model has been devised,
that is,
in a single transition multiple (up to three) vessels may depart from any number of fertile stars in  $\mathcal{S}$, ``concurrently'' settling multiple stars that will be added to $\mathcal{S}$ in order to generate the successor state $\mathcal{S}'$:

\begin{equation}
    \mathcal{S}' = \mathcal{S} \cup \bigcup_{(i \leq n_f)} \bigcup_{(j \leq 3-|\mathcal{O}_{k_i}|)} \left\{s_{o_j} \in \mathcal{N}_{k_i} \, |\, s_{k_i} \in \mathcal{S} \wedge |\mathcal{O}_{k_i}| < 3 \right\}
\end{equation}

The number of successors that can be generated in this way is much larger, as it depends on all possible combinations of departure and arrival stars, 
letting aside the travel time, which is assumed to be the minimum.
%
%
This number increases exponentially as the settlement tree grows; it also increases as the star density increases, that is,  at lower radii (less then 5 kpc).

Both memory and time constraints prohibit the generation of all possible successors.
%
%
 %
%
A constructive heuristic is proposed for generating up to $b_{f,max}$ successors from a single ancestor,
by using a well-defined set of rules for the selection of: a) fertile (departure) stars, b) the number of out-bound vessels from each departure star, and c) the neighbor (arrival) star for each vessel.

An iterative procedure has been devised to built the set of successors. 
Let $\mathcal{L}^i$ be a set containing the successors generated at the $i$-th iteration ($\mathcal{L}^0={\mathcal{S}}$ at the first iteration). 

\begin{enumerate}
    \item First, the fertile stars in the ancestor state $\mathcal{S}$ are ranked by using as ``fitness'' the number of missing stars (\emph{i.e.}, still to settle) in their own cell $(k_r, k_\theta)$ and in the 8 adjacent cells:
    \begin{equation}
        f = \sum_{p=-1}^{1} \sum_{q=-1}^{1} \left [ G(k_r + p, k_\theta + q) - H(k_r + p, k_\theta + q) \right ]
    \end{equation}
    Then, a fertile star $s_k$ is selected among the candidates (and removed from the list) by means of a biased sampling without replacement, proportional to $f$. In such a way, it is promoted the creation of branches directed towards the most empty cells of the zone.
    \item Now, $m_k\,|\mathcal{L}^i|$ new successors states are created,
    by appending $m_k$ differently generated set of offspring of star $s_k$ to each state in $\mathcal{L}^i$.
    The parameter $m_k$ is set according to the formula:
        \begin{equation}
    \begin{array}{cc}
         m_k &=
         \begin{cases}
         1.2\, \hat m & \mbox{if}\,\,\, f_k > \bar{f} \\
         \hat m & \mbox{otherwise}
         \end{cases}
    \end{array}
    \label{eq:nsucck}  
    \end{equation}
    with 
    $\hat m \in [\hat m_0, \hat m_0 + 4]$ an user-defined value, either constant or randomly sampled at any iteration, where the value $\hat m_0 = \sqrt[\leftroot{-3}\uproot{3}n_f]{b_{f,max}}$ ensures  a symmetrical successor tree; 
    $\bar{f}$, instead, is the average fitness of the fertile stars.
    The idea is to stress the exploration of empty cells by allowing a 20\% more 
    successors states to be generated from stars with a fitness larger than the average.

    \item The number of offspring $\beta_k$ of the star $s_k$ to insert in each new successor is randomly sampled in~$\left [1, 3 - | \mathcal{O}_k | \right]$, according to an uneven probability distribution that aims at increasing the number of settled stars as soon as possible.
    \item 
    Finally, each offspring star $s_o$ must be selected within the neighborhood $\mathcal{N}_k$ of the parent star $s_k$.
    Each offspring star $s_o$ is chosen either as the closest in time to the parent star $s_k$, or as the one that realizes the maximum displacement in $\theta_f$ per unit time.
    Randomly alternating these two rules allows to combine 
    large, but time consuming, horizontal movements, and
   ``quick reproduction'' phases corresponding, in most cases, to vertical movements
    in the ($r$, $\theta_f$) plane. 
    As a result, the following \textit{throw and multiply} rule for offspring selection is defined:
    \begin{equation}
    \begin{array}{cc}
         s_o &=
         \begin{cases}
         \argmin_{s_h \in \mathcal{N}_k} (t_h - t_k) & \mbox{if}\,\,\, p_r < 0.5 \\
         \argmax_{s_h \in \mathcal{N}_k} \frac{|\Delta \theta_f|}{t_h - t_k} & \mbox{otherwise}
         \end{cases}
    \end{array}
    \label{eq:succ_choice}  
    \end{equation}
    being $p_r$ a random number between 0 and 1 and $\Delta \theta_f$ the difference in final polar angle between star $s_h$ and $s_k$.
    If the stars closest in time are more than one, one of them is randomly picked when $p_r < 0.5$.
    At this point, $s_o$ is removed from $\mathcal{N}_k$ and another offspring star is selected with the same rule, until $\beta_k$ offspring are generated.
    The $\beta_k$ offspring are appended to one of the $m_k\,|\mathcal{L}^i|$ new successors under construction, and the points 3. and 4. are repeated for each of them.
\end{enumerate}

The procedure is repeated until: i) it is no longer possible to create new successors (that is, all neighbourhoods are empty or there is no available fertile star) or ii) a number of successors equal to $b_{f,max}$ has been generated.
In both cases, all the successors in the last evaluated iteration, $\mathcal{L}^{last}$, are retained.  

\begin{figure} [h!]
\captionsetup[subfigure]{justification=centering}
    \centering
    \begin{subfigure}{0.49\textwidth}
    \centering
    \includegraphics[width = 0.8\textwidth]{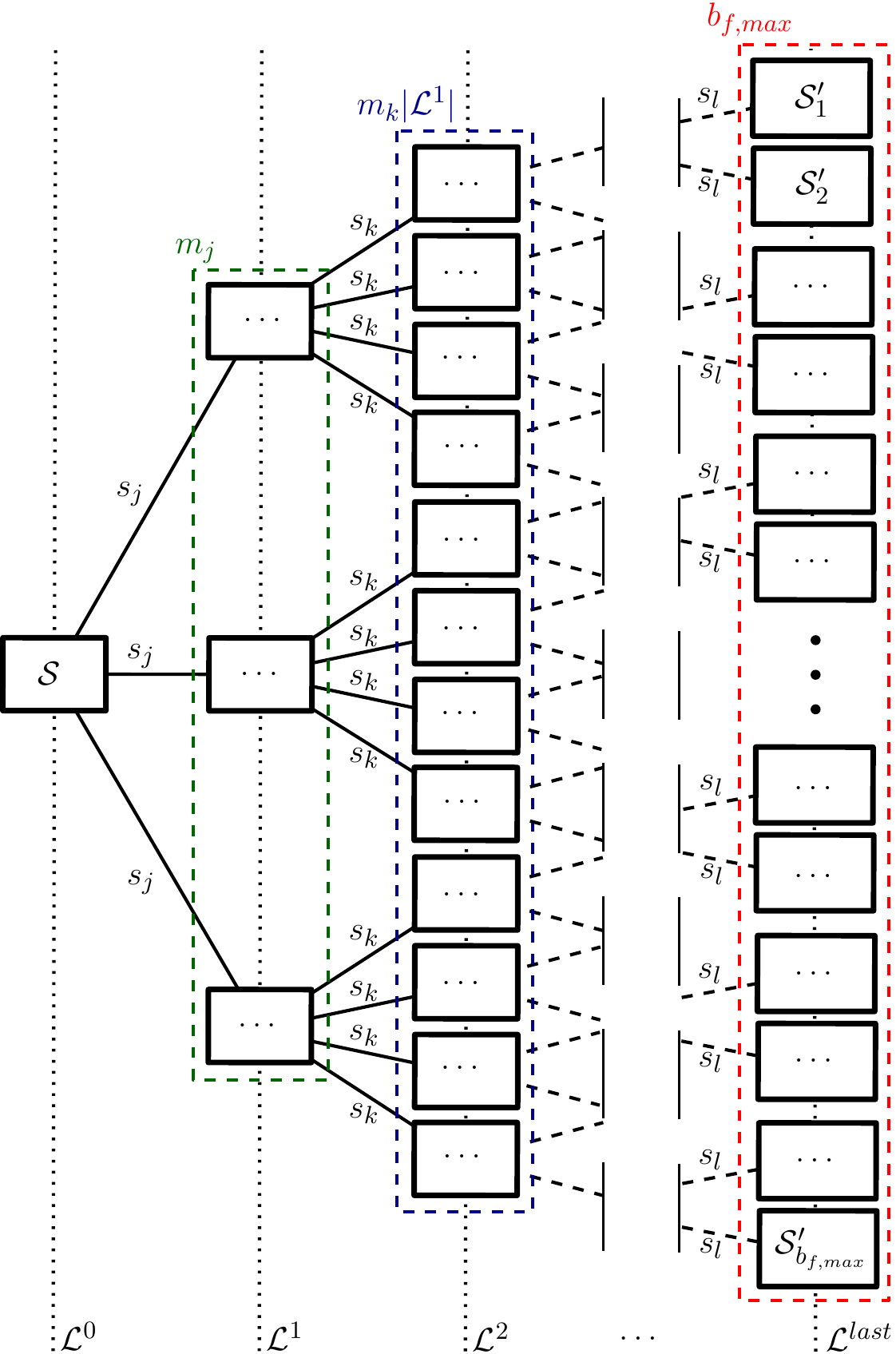}
      \caption{Successor tree. The star $s$ reported on each branch identifies the fertile star selected for offspring generation in the corresponding iteration.}
   \label{fig:transition_tree}
    \end{subfigure}
    \begin{subfigure}{0.49\textwidth}
    \centering
    \includegraphics[width = 0.93\textwidth]{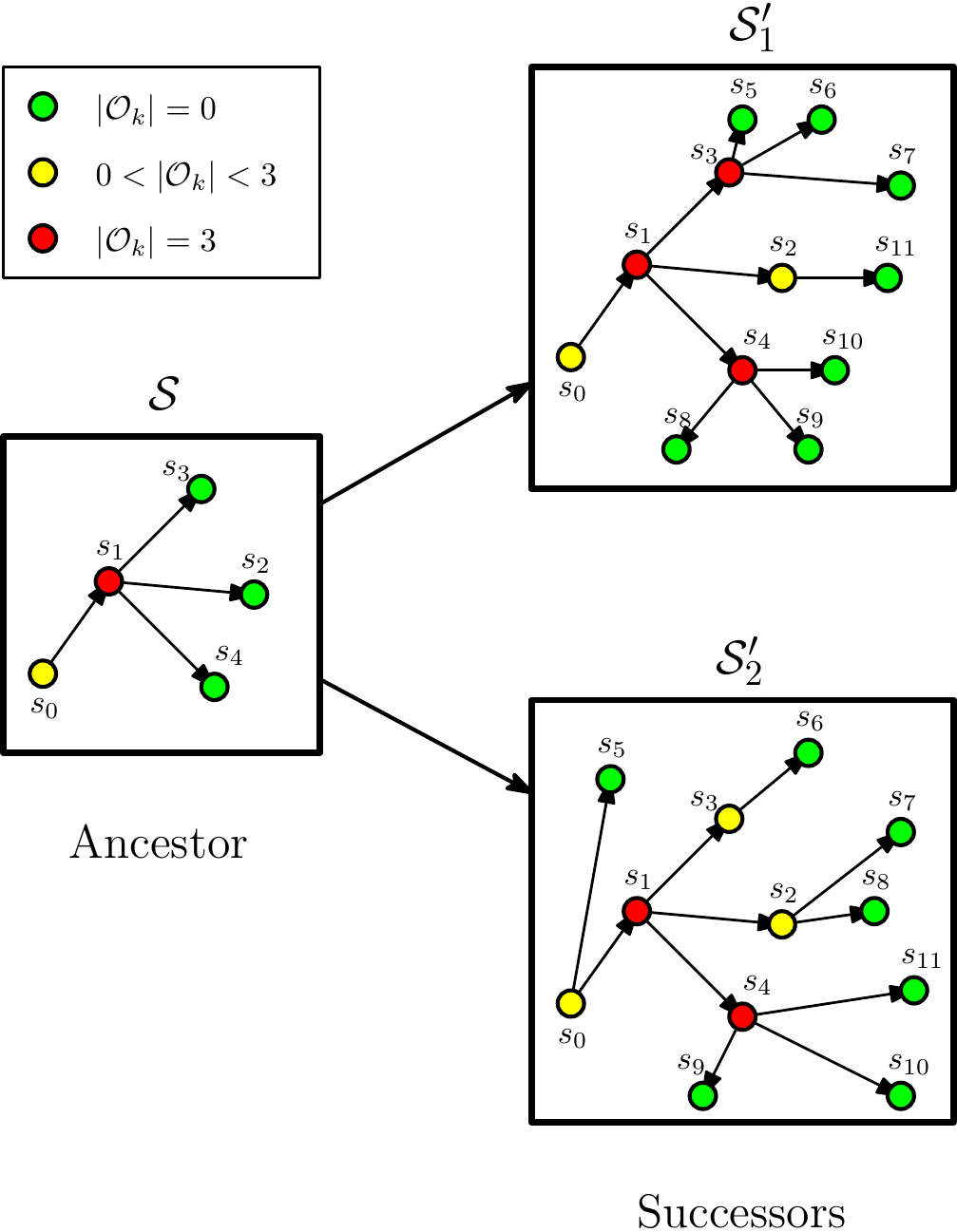}
      \caption{Ancestor state $\mathcal{S}$ and two successors $\mathcal{S}'_1$ and  $\mathcal{S}'_2$.}
    \label{fig:transition_state}
    \end{subfigure}
    \caption{Transition model: successor tree (a), ancestor and successor states (b).}
    \label{fig:transition}
\end{figure}

An example of the successor tree generated in a single transition from the ancestor state $\mathcal{S}$
is shown in Figure~\ref{fig:transition_tree}; the settlement trees corresponding to the ancestor state $\mathcal{S}$ and two of its successors, $\mathcal{S}'_1$ and  $\mathcal{S}'_2$, are shown in Figure~\ref{fig:transition_state}.

\subsection{Beam Best-First Search} \label{sez:BBFS}
A multi-settler stochastic Beam Best-First Search (BBFS) has been developed for solving the unusual search problem defined in the previous section,
as a generalization of the (single-spacecraft) ``Stochastic Beam'' algorithm proposed by Sim{\~o}es et al. \cite{simoes2017multi} for  handling the multi-rendezvous problem released in GTOC 5. 
Solution methodologies based on the beam search have been demonstrated to be successful also in other GTOCs, characterized by a wide search space \cite{grigoriev2013choosing, izzo2013search, izzo2016designing, petropoulos2017gtoc9}.

Being the beam search a \textit{tree-search} algorithm, it tries to approach the goal gradually building a so-called \textit{search tree}, with the initial state as root: 
the \textit{nodes} of the tree correspond to the explored states in the state space of the problem, while the \textit{branches} correspond to state-to-state transitions.

At each iteration, the algorithm picks the first \textit{leaf} state in the \textit{frontier}, that is, the list of all the generated, but yet unexplored, nodes at any given point of the search, and \textit{expands} it by generating a set of \textit{successor} states. 
In the proposed beam search variant, termed \textit{best-first}, the frontier is modeled as a priority queue ordered by $\phi$. So, at any iteration, the best state is used for generating $b_{f,max}$ successors states. A \textit{probabilistic branching} is now exploited so as to avoid to fill the frontier only with the newly generated nodes:
with a probability $p$, the best $b_f$ successor states, according to $\phi$, are retained and appended to the frontier; otherwise, the $b_f$ successors are chosen by means of a biased sampling without replacement proportional to $\phi$.
Subsequently, the frontier, now containing a (possibly large) number of states evaluated at different iterations, is sorted according to the cost function $\phi$ and the best $b_w$ are kept in memory, while the others are discarded.
The now-exhausted ancestor state is removed from the search-tree frontier, and a new iteration may begin.
The algorithm terminates as soon as successors can no longer be generated from any of the states in the frontier, or a maximum, pre-defined, number of nodes have been explored. The best-so-far state is eventually elected as problem solution. 

\begin{figure} [h!]
\centering
\begin{subfigure}{0.48\textwidth}
  \centering
  \includegraphics[width=\textwidth]{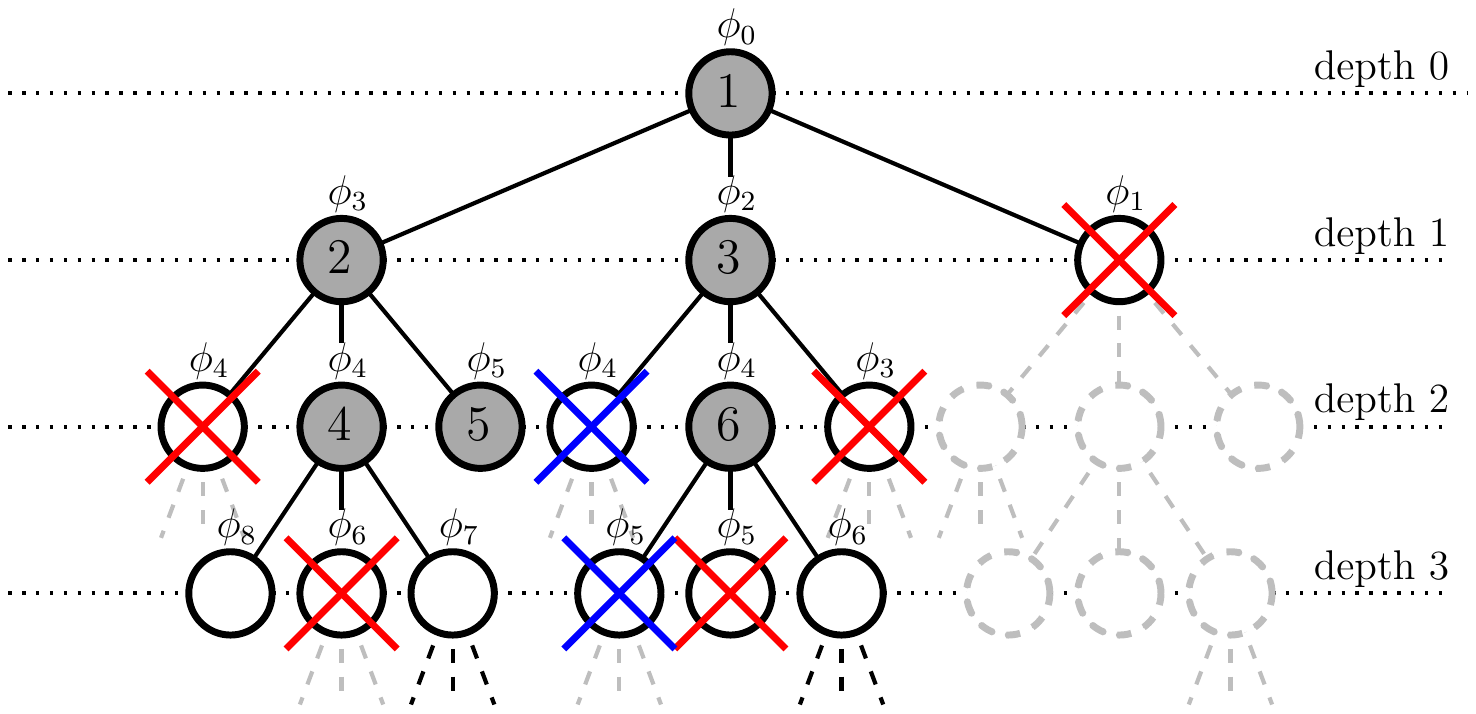}
  \caption{Beam (breadth-first) Search.}
  \label{fig:BS}
\end{subfigure}
\begin{subfigure}{0.48\textwidth}
  \centering
  \includegraphics[width=\textwidth]{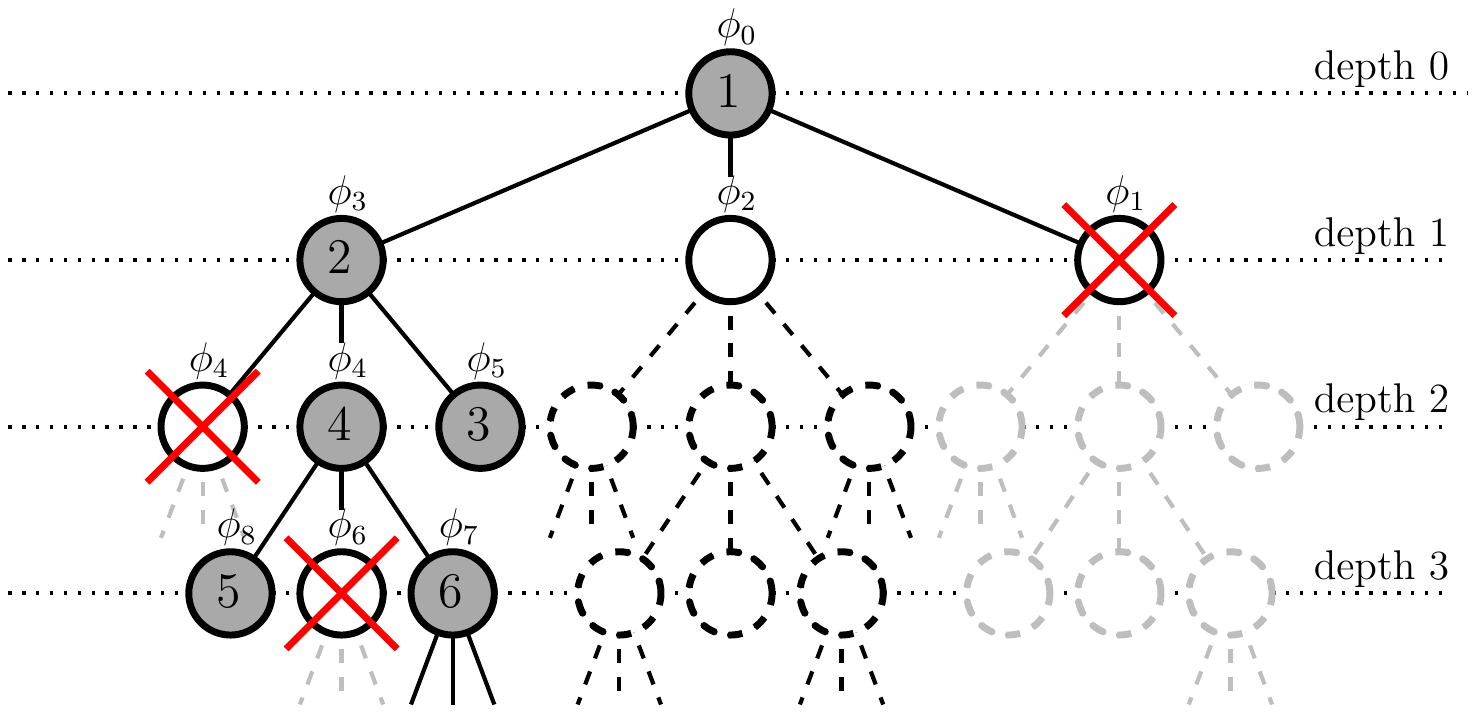}
  \caption{Beam Best-First Search.}
  \label{fig:BBFS}
\end{subfigure}
\caption{Comparison between six iterations of a breadth-first and best-first version of the Beam Search on the same problem, by supposing $b_f = 2$ and $b_w = 3$. In particular, grey-filled nodes have been explored (with the number indicating the order); white-filled nodes have been generated but not yet explored; black-dashed nodes are yet to be generated; grey-dashed nodes are the successors of a (crossed) node pruned by the beam procedure. More precisely, red-crossed nodes are excluded by the probabilistic branching ($b_f$), while blue-crossed nodes by the frontier re-sizing ($b_w$). The values of the cost function $\phi$ are labeled in decreasing order: $\phi_j < \phi_i$ if $j > i$.}
\label{fig:tree_search}
\end{figure}

The best-first (or uniform-cost \cite{russell2016artificial}) variant has been here adopted, in place of the ``standard'' breadth-first search, in order to provide good-quality answers to each set-covering problem as soon as possible: indeed, as shown in Figure~\ref{fig:tree_search}, the best-first variant is able to cross freely the different levels of the tree and quicker generate ``deeper'' nodes (\emph{i.e.}, states with a larger number of settled stars), that we expect have a better solution score. 
In this way, settlement trees of up to 100 stars were generated in a few minutes on a standard laptop.

The hyper-parameters of the BBFS procedure were tuned ``by hand'' during the competition, by trying different values and comparing the correspondent results. At the end, the following values resulted to be the best ones for the great part of the analyzed zones: $b_w = 20000$, $b_{f,max} = 20000$, $b_f = 1000$.
In order to mitigate potential issues due to the ``non-exhaustiveness'' and stochastic nature of the BBFS, several runs were performed in parallel, and the best-found solution saved.

The proposed procedure has been repeated over all zones,
iteratively changing the target occupancy distribution according to the attained results,
in order to have an overall number of settled stars in each row and column as close as possible to the one stated by the global, 1020-star, target distribution devised beforehand.
Eventually, a settlement tree with $1013$ stars was found at the end of this process.

\section{Ideas for solution improvement}

A number of refinement strategies were investigated during the competition in order to improve the BBFS solution, 
as soon as a 1020-star settlement tree were available.
The three most effective are here presented.

\subsection{Concurrent optimization of transfer times}

Once the layout of the settlement tree has been defined, a large score boost can be obtained by a proper 
optimization of the intermediate departure and arrival times as well as  by moving to 90 Myr the settle epoch of all terminal stars (i.e., stars without offspring).

A rudimentary sequential quadratic programming optimization has been  carried out.
Let $\delta = \delta(j)$ be the additional transfer time admitted for the $j$-th transfer leg composing a settlement tree of $N$ legs (no special sort of the transfers is required).
Modified departure and arrival times $t_f(j)$, $t_0(j)$ are evaluated as
\begin{align}
    	 t_f(j) &= \bar t_f(j) + \delta(j) \\
	 t_0(j) &= \bar t_0(j) + \sum_{i=1}^{N}{ M(i,j) \delta(i)} 
\end{align}
where over-lined variables refers to the unperturbed quantities,
and $M(i,j)\in\left\{0,1\right\}$ is an $N\times N$ matrix indicating if the arrival star of the $i$-th transfer arc
is the same as the departure star of the $j$-th arc ($M(i,j) = 1$) or not ($M(i,j) = 0$).
%
A quadratic approximation of the total used $\Delta V$ is built using the derivatives of each transfer cost $\Delta V(j)$ with respect to its departure and arrival time: 
\begin{align}
    \begin{split}
	&  \Delta V(j) =  \overline{ \Delta V}(j) 	+ \frac{\partial\Delta V}{\partial t_0} \left(t_0(j) - \bar t_0(j) \right)
	+\frac{1}{2} \frac{\partial^2\Delta V}{\partial t_0^2} \left(t_0(j) - \bar t_0(j) \right)^2  \\
	& \phantom{\Delta V =   \Delta V }	 
	+ \frac{\partial\Delta V}{\partial t_f} \left(t_f(j) - \bar t_f(j) \right)
	+\frac{1}{2} \frac{\partial^2\Delta V}{\partial t_f^2} \left( t_f(j) - \bar t_f(j) \right)^2   
	\end{split}
\end{align}

A (local) minimization problem is then formulated as:
\begin{align}
 \min_{\delta_{min} \le \delta  \le \delta_{max}} 	& { \sum_{j=1}^{N} { \Delta V(j) }    }    \\
	  \mbox{s.t.}\,\,\,\, &    \Delta V(j) < \Delta V_{tot}^{max}  \quad \forall \, j \in \left\{1,\ldots, N\right\}
\end{align}
where $\delta_{min}$, $\delta_{max}$ are sufficiently close bounds that ensure the validity of the quadratic approximation;
$\delta_{min} = -0.2$ Myr and  $\delta_{max} = 1$ Myr 
have been selected according to the author experience.
This procedure is repeated, each time updating the reference value with the newly obtained ones.
A few iterations 
are usually sufficient to acquire a significant reduction of the  overall mission $\Delta V$.

%

\subsection{Explosion}

The overall settlement tree, obtained as the sum of the solutions of the different set covering problems,  may still grow in size (\emph{i.e.}, in number of stars) without hindering the 
distribution score $J_2$.
In fact, with respect to the previously discussed set covering problem, 
this settlement tree already covers all the cells of interest.
Therefore, it is now possible to work on a simpler cell-by-cell basis, while keeping the distribution as uniform as desired.

The explosion algorithm is here presented.
%
Let $\mathcal{S} = \left\{ s_1 , s_2, \dots , s_N \right\}$ be the set of settled stars so far. Each star $s$ is characterized by a certain number of available Settler Ships $n_A(s) = 3 - |\mathcal{O}(s)|$ that can still depart and reach new stars.
Let $x( k_r , k_\theta )$ be an integer variable indicating the number of additional stars to settle in cell $( k_r , k_\theta )$.
This variable is non-negative and upper bounded by the maximum number of stars $x_{\text{max}}( k_r , k_\theta )$ that can be settled in the cell.
Working on a cell-by-cell basis, that is, only transfers within a cell are allowed,  
a rough estimate of this upper bound is given by the sum of the available settlers in that cell:
\begin{equation}
	x_{\text{max}}( k_r , k_\theta ) = \sum\limits_{s \in ( k_r , k_\theta )} n_{A}(s)
\end{equation}

If newly settled stars are reached early in time
a new set of Settler Ships is available for further colonization of the galaxy, thus generating an entire new ``army''.
%
%
 %
Therefore, a more accurate estimate of the maximum number of new stars per cell,
given the problem time/impulse magnitude constraints,
is given by:
\begin{equation}
	x_{\text{max}}( k_r , k_\theta ) = \sum\limits_{s \in ( k_r , k_\theta )} \left( n_{A}(s) \sum\limits_{k = 0}^{n_{\text{gen}}} 3^{k} \right)
\end{equation}
where $n_{\text{gen}}$ is the number of following generations. 
It can be computed as:
\begin{equation}
	n_{\text{gen}} = \left\lfloor \frac{t_{end} - t(s)}{\Delta t + t_w} \right\rfloor
\end{equation}
where $t(s)$ denotes the settling time of the departure star $s$, $t_w$ is the waiting time before a Settler Ship can depart from a settled star, $t_{end}$ denotes the mission final time and $\Delta t$ is the average transfer time (assumed equal to 6 Myr). 

The (putative) optimal number of additional stars $x(k_r, k_\theta)$ to settle in each cell is found as the solution of an Integer Linear Programming (ILP) problem
that aims at minimizing the row/column error on the novel star distribution  with respect to the  desired target one, given
i) the current star distribution $N_0(k_r, k_\theta)$, 
ii) a target row/column star distribution $N^{des}_r(k_r)$,\,$N^{des}_\theta(k_\theta)$, and
iii) a rough estimate of the maximum number of stars the could be settled in a given cell $(k_r, k_{\theta})$ by working on a cell-by-cell basis.
The problem can be stated as:
\begin{align}
	&\min_{0<X<X_{max}} \sum\limits_{k_r = 1}^{30} \Delta_r(k_r) + \sum\limits_{k_\theta = 1}^{32} \Delta_\theta(k_\theta)  
\end{align}
where $\Delta_r(k_r)$, $\Delta_\theta(k_\theta)$ are auxiliary error functions defined as the difference between the settled stars (both before and after the explosion) and the desired ones in a ring or, respectively, in a slice:
\begin{align}
\Delta_r(k_r) &= \left\lvert \sum_{k_\theta = 1}^{32}{\left( x\left( k_r , k_\theta \right) +  N_0\left( k_r , k_\theta \right)\right)} - N_r^{\text{des}}( k_r )  \right\rvert \\
\Delta_\theta(k_\theta) &= \left\lvert \sum\limits_{k_r = 1}^{30}{  \left( x( k_r , k_\theta ) + N_0( k_r , k_\theta )\right)} -  N_\theta^{\text{des}}( k_\theta ) \right\rvert
\end{align}

By repeatedly solving for a target distribution that is a multiple of the theoretical optimal 512-star target distribution, one hopes to find a refined target distribution that features a much larger number of stars without breaking the initial regularity.
Once the putative optimal number of new stars to settle in each cell has been obtained,
a greedy search is carried out on a cell-by-cell basis in order to  
find the actual settlement tree that realizes that distribution.
All possible combinations of minimum-time transfers between already settled stars and available stars in that cell are considered.
Fertile stars are sorted by their settle times; those with lower settle time are used first in order to generate as soon as possible new vessels that could possibly reach other stars before the end of the mission time.

\subsection{Pruning}

In most cases, the attained settlement tree is not as uniform as expected.
This may happen for a variety of reasons. As an example, 
the BBFS 
may have favoured a group of stars very close each other, 
or the explosive algorithm my have failed in uniformly expanding the original settlement tree.
In either case, a  non-uniform star distribution could severely hinder the overall score.
Therefore, in order to re-gain a more uniform star allocation,
the settlement tree undergoes a (greedy) removal process, \emph{i.e.}, a ``pruning":
for each star,  the increase in the merit index $J_2$  associated with its elimination  is computed;
the star leading to the largest score improvement is actually removed.
This procedure is iteratively repeated as long as the score function keeps improving.
This step should always be followed by a concurrent optimization of all transfers time, as reducing the number of stars may also reduce the overall propellant consumption.

\section{Numerical Results}

\begin{figure}[h!]
\centering
\captionsetup{justification=centering}
\begin{minipage}{0.45\textwidth}
    \centering
    \includegraphics[width = \textwidth]{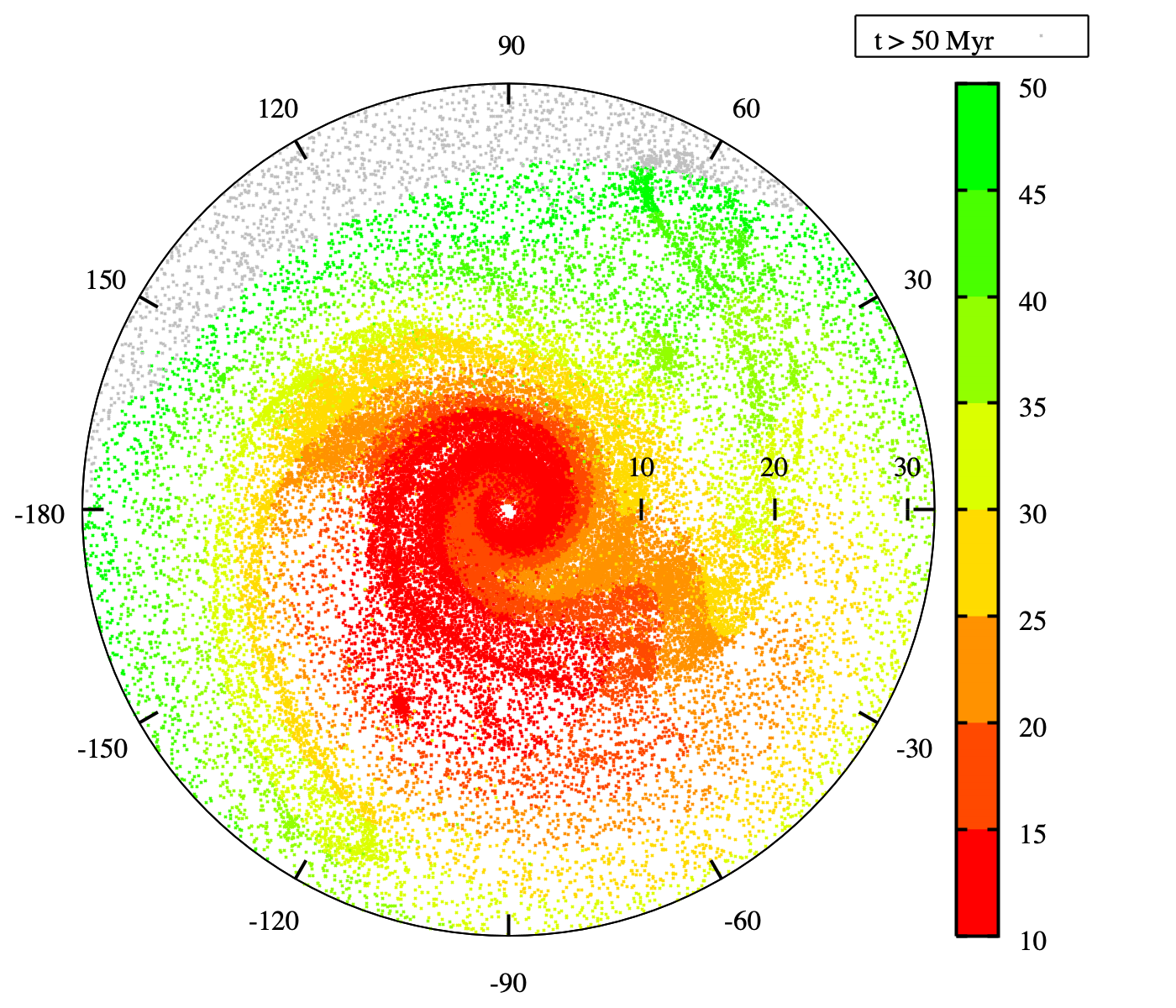}
    \caption{Reachable stars in the galaxy as a function of the transfer time for the Fast Ships. }
    \label{fig:fast}
\end{minipage}
\begin{minipage}{0.54\textwidth}
    \centering
    \includegraphics[width = \textwidth]{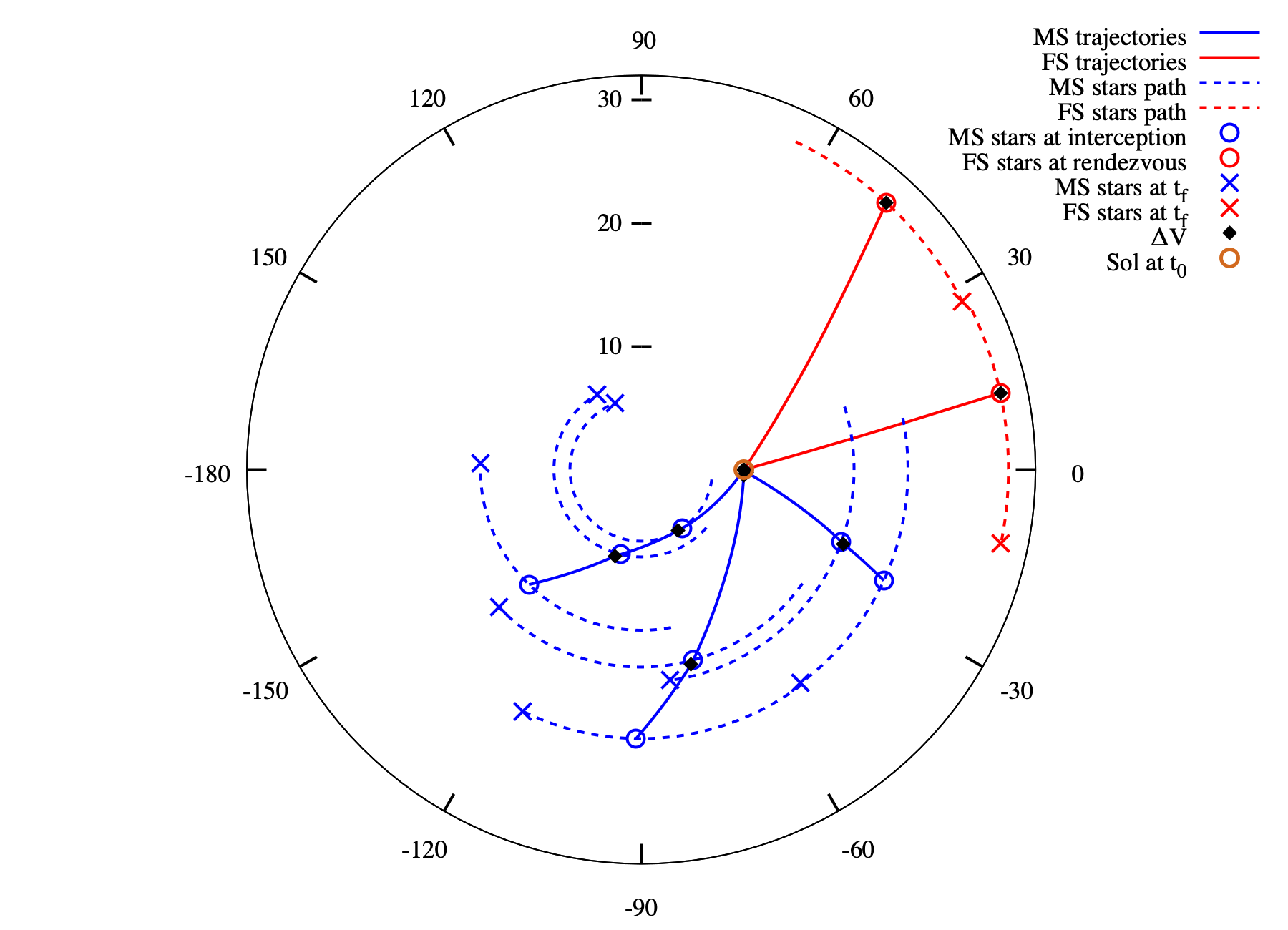}
    \caption{Trajectories of Mother Ships and Fast Ships: view from Galactic North. }
    \label{fig:MS_FS}
\end{minipage}
\end{figure}

\begin{figure}[!h]
\centering
\begin{subfigure}{0.49\textwidth}
  \centering
    \includegraphics[width=1.1\textwidth]{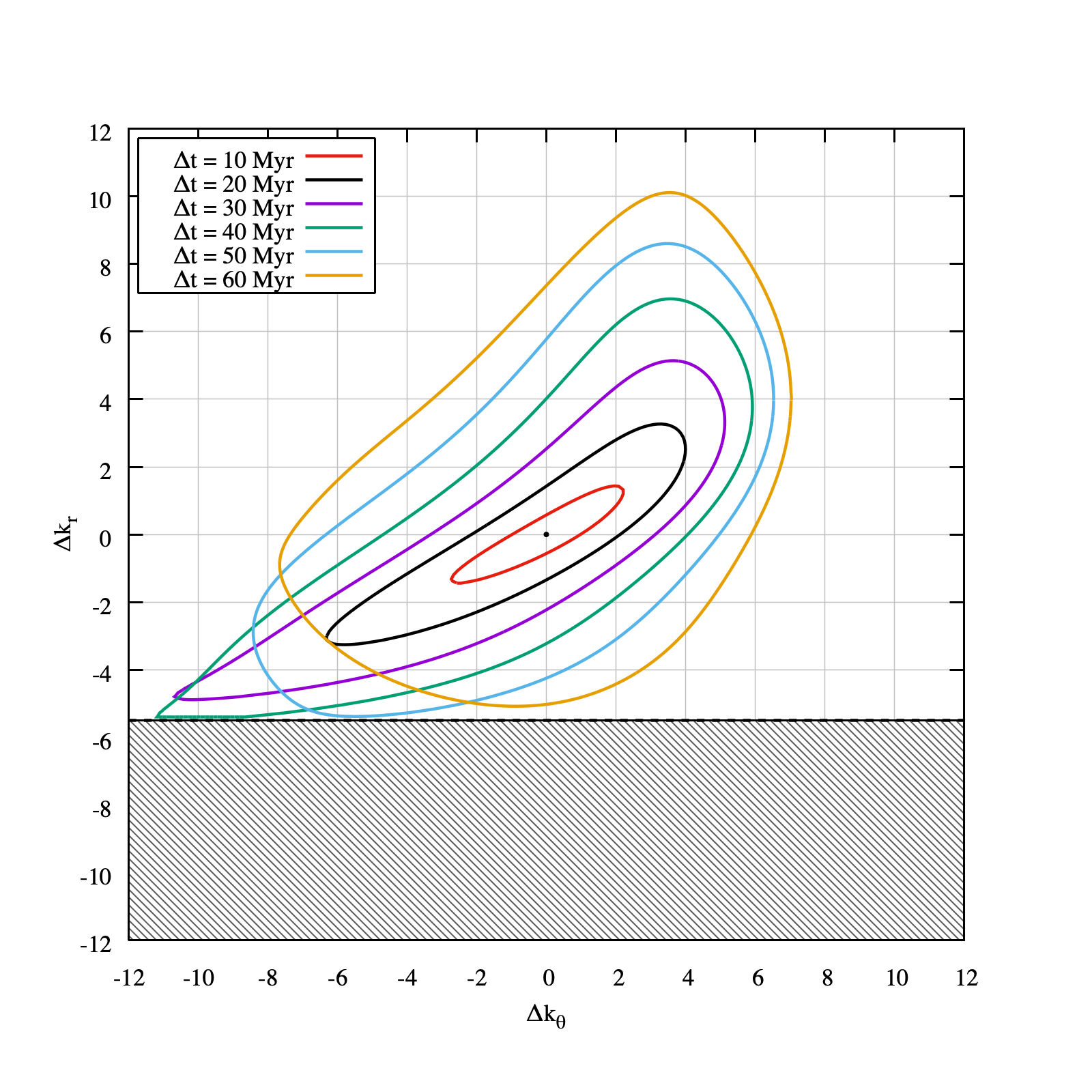}
    \caption{$r_0 = 7.5\,kpc$.}
    \label{fig:Map19}
\end{subfigure}
\begin{subfigure}{0.49\textwidth}
  \centering
    \includegraphics[width=1.1\textwidth]{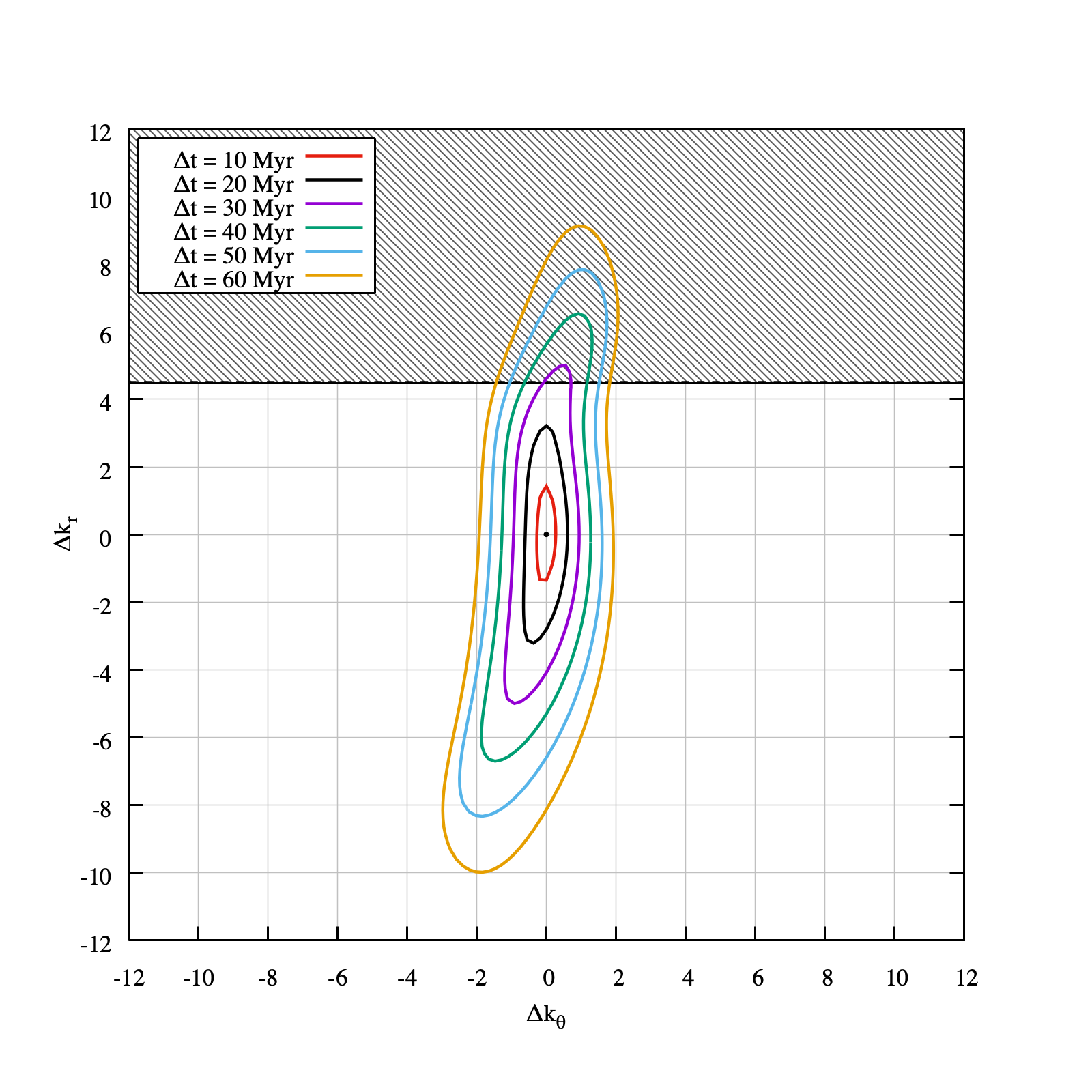}
    \caption{$r_0 = 27.5\,kpc$.}
    \label{fig:Map25}
\end{subfigure}
\caption{Mobility maps.
Contour lines are represented
for a vessel departing from two different radii at $t_0 = 30$ Myr, as a function of the transfer time $\Delta t$. Patterned zones indicate regions beyond the given boundaries of the galaxy.}
\label{fig:mobilityMap}.
\end{figure}

Some numerical results that document the adopted procedure are here presented.
%
%
At the beginning of the competition a large effort has been devoted to the analysis of the mobility of Fast Ships and Settler Ships. 
Figure~\ref{fig:fast} highlights the stars that is possible to reach with a Fast Ship as a function of the allowed travel time. %
Mobility charts as those shown in Figure~\ref{fig:mobilityMap} 
were instead used for understanding the capability of a  Settler Ship to reach nearby stars, with respect to the allowed travel time.
From this analysis, it was apparent
the scarce lateral mobility at higher radii, 
somehow opposed by the great lateral mobility at lower radii, and the usually good radial mobility at any radii.
This analysis was also fundamental for the definition of the 
\textit{throw and multiply} strategy exploited in the BBFS procedure.

\begin{figure} [h!]
\centering
\begin{subfigure}{0.49\textwidth}
  \centering
    \includegraphics[width=\textwidth]{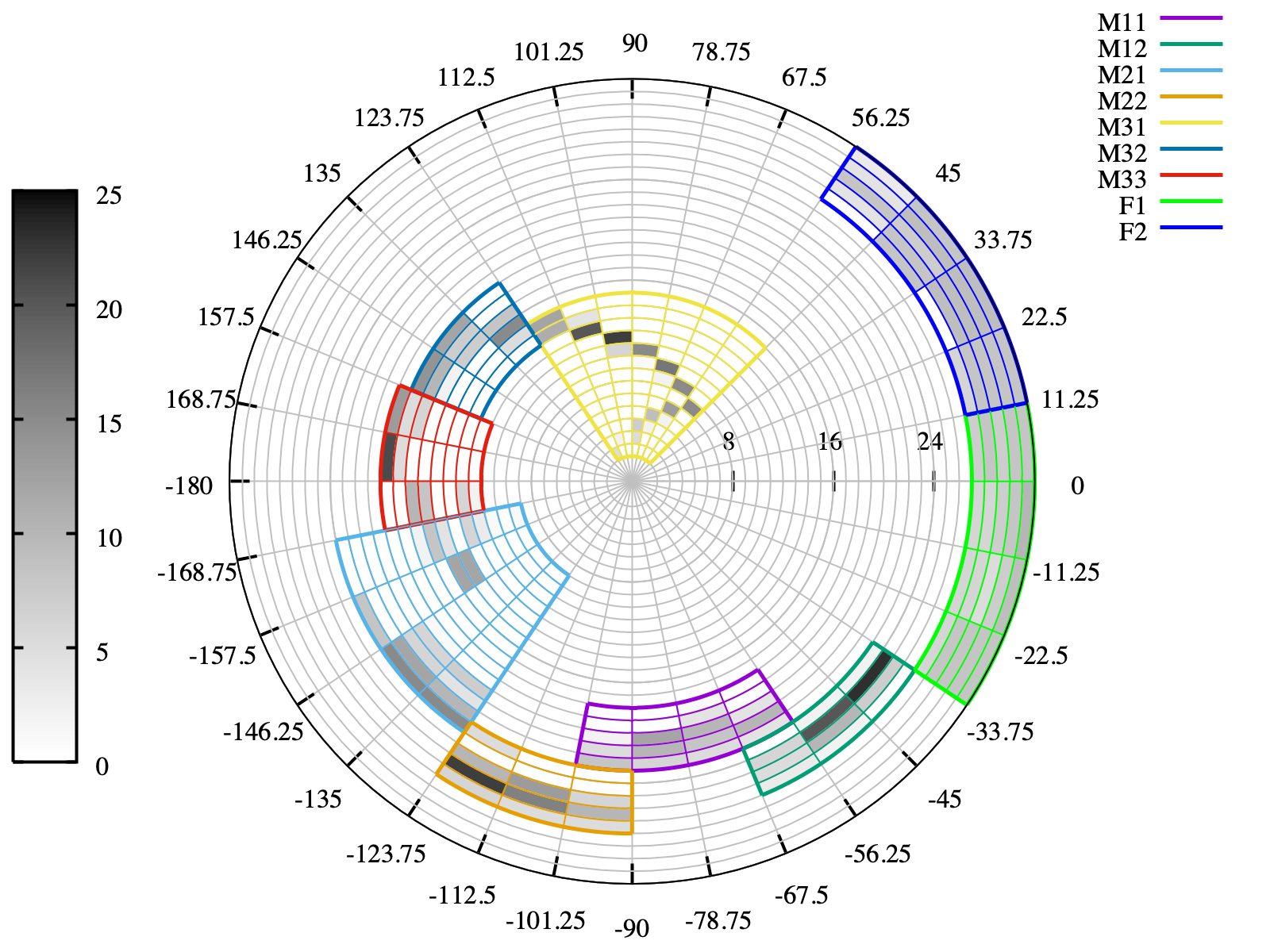}
    \caption{Target distribution in the ($r$, $\theta_f$) plane.}
    \label{fig:zonesTar}
\end{subfigure}
\begin{subfigure}{0.49\textwidth}
  \centering
    \includegraphics[width=\textwidth]{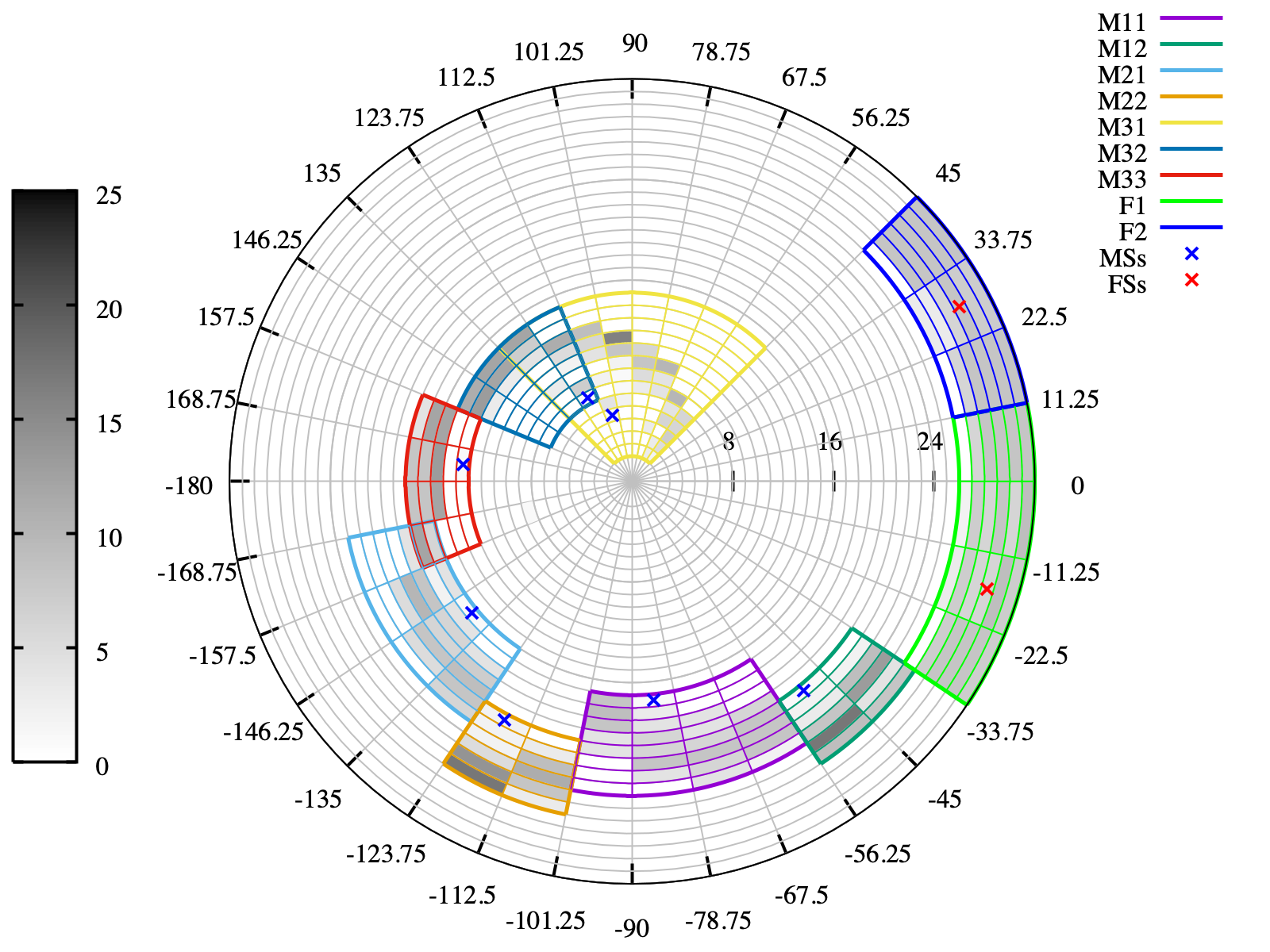}
    \caption{Obtained distribution in the ($r$, $\theta_f$) plane. }
    \label{fig:zonesFin}
\end{subfigure}
\caption{Identified zones and cells occupation: target distribution (a), obtained distribution (b). The grey-scale refers to the number of stars to settle in any cell. Blue and red crosses identify the cells where the ``root'' stars encountered by Mother Ships and Fast Ships, respectively, are located in the final solution.}
\label{fig:zones}
\end{figure}

After a long reasoning, 
nine (almost) non-overlapping zones in the plane ($r$, $\theta_f$) were figured out.
Each zone refers to a ``root'' star either provided by a Fast Ship or by a Mother Ship: one of the Mother Ships intercepts three stars and the others two stars each; their trajectories are reported in Figure~\ref{fig:MS_FS}.
Independent settlement trees departed from each root star.
These  zones are listed in Table~\ref{tab:zonesTab}, and shown in Figure~\ref{fig:zonesFin} in the ($r$, $\theta_f$) plane. 
The following nomenclature has been used to label each zone:
the initial letter indicates if the star is encountered by a Mother Ship ($M$) or a Fast Ship ($F$), respectively; the first number uniquely identifies the ship, and, the second number (only for Mother Ships) corresponds to the settled star number (ordered by settling time).

Figure~\ref{fig:zonesTar} provides a visual representation of the initial (tentative) 1020-star distribution on the ($r$, $\theta_f$) grid.
As more information about the problem were gathered, 
the target occupancy distribution was updated, up to the very last day, bringing to the final distribution shown in Figure~\ref{fig:zonesFin}.

\begin{table}[h!]
    \centering
    \begin{tabular}{c c c c c c c c c c}
\hline\hline
	 $Z$        &$M11$&	$M12$&	$M21$& 	$M22$& 	$M31$& 	$M32$& 	$M33$& 	$F1$ & 	$F2$ \\ 
\hline
$k_r$	    & 8 $\div$ 11 &	12 $\div$ 13&	2 $\div$ 5	 &  6 $\div$ 7	 &  21 $\div$ 28&	27 $\div$ 30&	31 $\div$ 2 &	14 $\div$ 17&	18 $\div$ 20\\
$k_\theta$	&16 $\div$ 23&	20 $\div$ 25&	15 $\div$ 21&	20 $\div$ 25&	1 $\div$ 13 &	6 $\div$ 13 &	12 $\div$ 16&	25 $\div$ 30&	25 $\div$ 30\\
\hline\hline
   \end{tabular}
    \caption{Zone names and ranges.}
    \label{tab:zonesTab}
\end{table}

\begin{figure} [h!]
\captionsetup[subfigure]{justification=centering}
    \begin{subfigure}{0.49\textwidth}
    \centering
    \includegraphics[width=\textwidth]{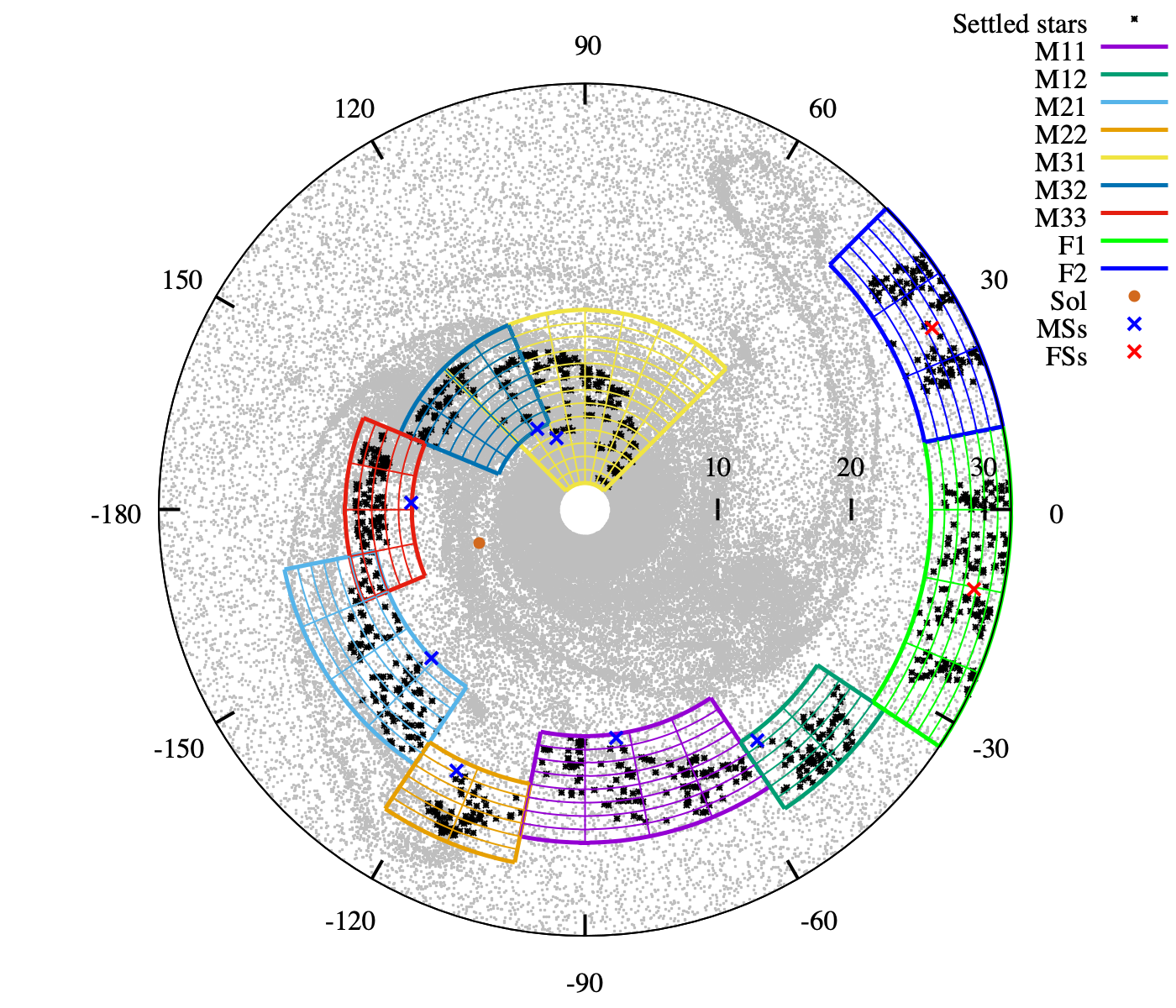}
    \caption{Submitted solution: view in the ($r$, $\theta_f$) plane.}
    \label{fig:SolGrid}
    \end{subfigure}
     \begin{subfigure}{0.49\textwidth}
    \centering
    \includegraphics[width=\textwidth]{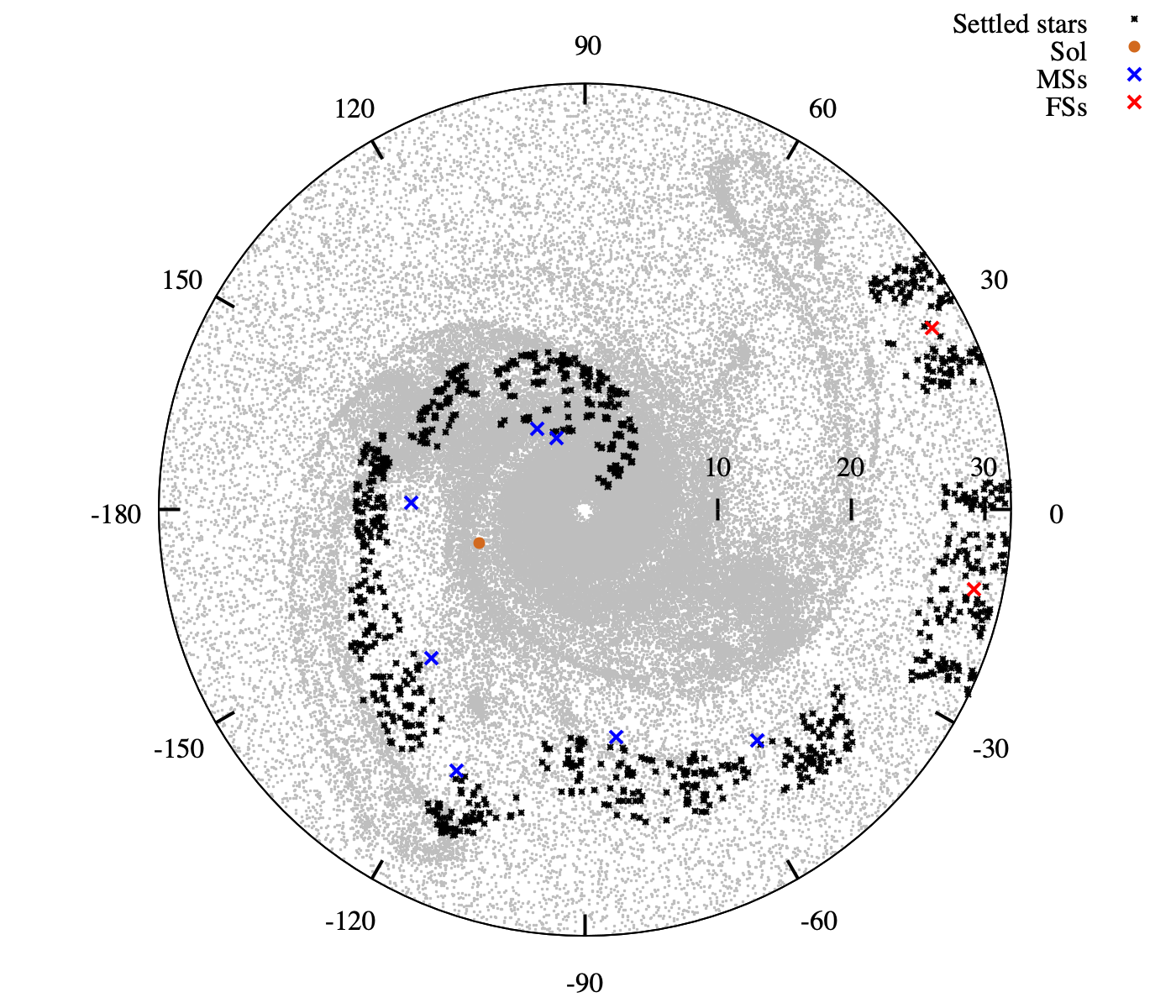}
    \caption{Submitted solution: view from Galactic North.}
    \label{fig:Enhancedsol}
    \end{subfigure}
    \caption{Final submitted settlement tree. The nine search zones are highlighted.}
    \label{fig:solution}
\end{figure}

At the end, the final submitted solution, shown in Figure~\ref{fig:solution}, is characterized by 1013 stars, ``uniformity'' merit index $J_2 =  630.069$ and ``propulsive'' merit index $J_3 = 1.504$, resulting in an overall score $J = 946.451$ and in the \nth{7} place in the leaderboard.

\subsection{Enhanced solution}

\begin{figure} [h!]
    \centering
    \includegraphics[width = \textwidth]{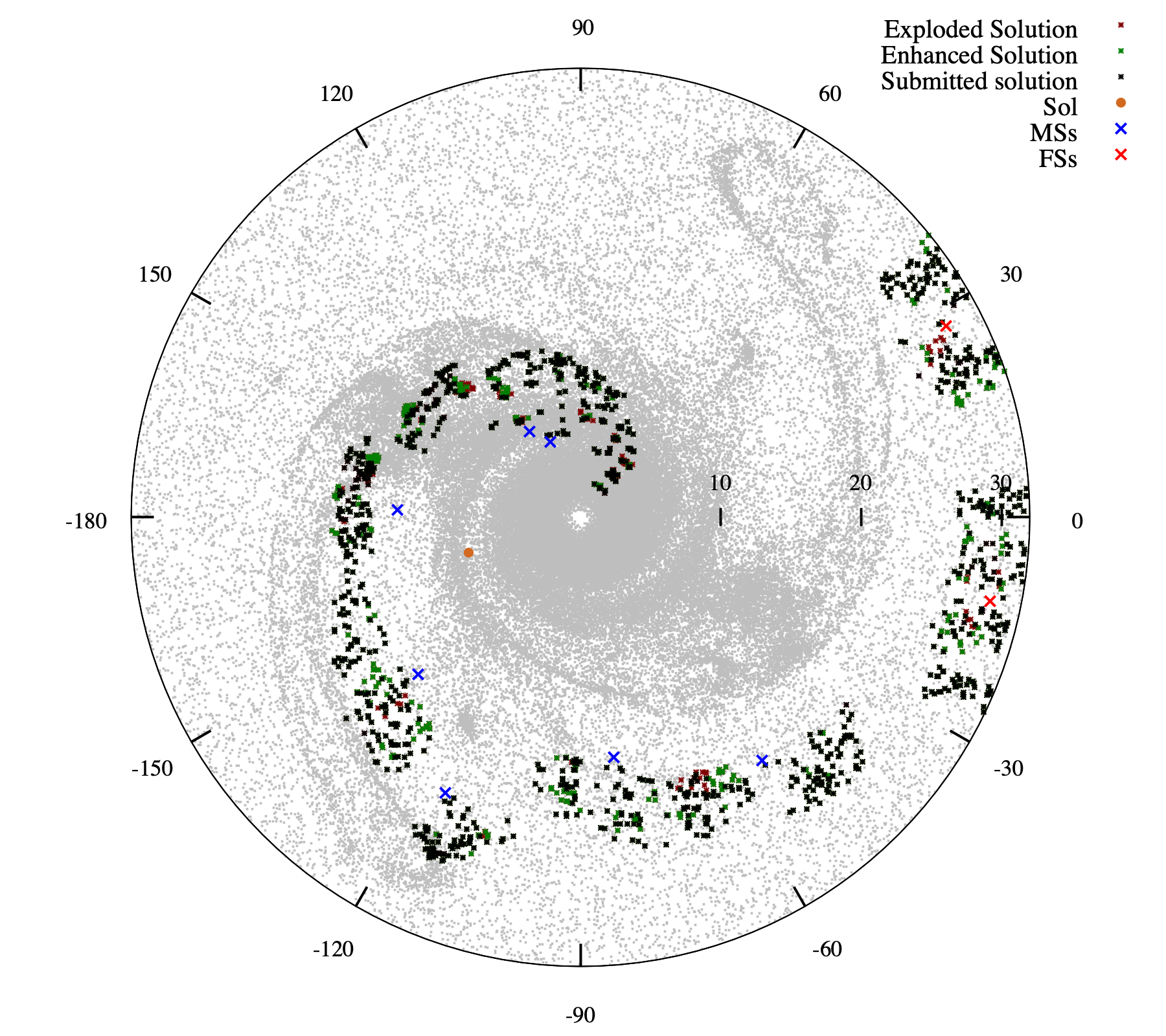}
    \caption{Comparison between the obtained solutions.}
    \label{fig:sol_complete}
\end{figure}

Due to the competition limited time, 
the submitted solution did not take advantage of the presented explosion and pruning techniques.
If these two ideas are applied to the same settlement tree, 
by first increasing  as most as possible the number of settled stars (explosion)  and then removing unnecessary settlements until $J_2$ reaches a local optimum (pruning), 
one attains a significantly better solution.
Figure~\ref{fig:sol_complete} presents 
the starting 1013-star settlement tree,
the ``exploded" solution with 1377 stars,
and the ``pruned'', final, solution with 1220 stars,
with merit index $J_2 = 772.761$ and $J_3 = 1.5531$,
that is an overall score $J = 1200.145 $. 
This result suggests that the proposed procedure is able to attain results competitive with the other top teams.

\begin{table}[h!]
    \centering
    \resizebox{\textwidth}{!}{
    \begin{tabular}{l c c c c c c c r}
    \hline\hline
         Solution          & $N$   & $E_r$ & $E_\theta$ & $J_2$ & $\Delta V_{used}$ [$10^3$ km/s]  & $\Delta V_{max}$ [$10^3$ km/s]  & $J_3$  & $J=J_2*J_3$ \\
    \hline 
    Submitted\phantom{$^\star$}  & 1013   & 1.88968 & 4.11965& 630.069 & 271.42  & 408.20 & 1.504  &   946.451 \\
    Exploded$^\star$  & 1375  & 3.67621 & 3.50581 & 691.814 &  405.65& 550.50 &  1.3571 &  938.839\\
    Enhanced\phantom{$^\star$}  & 1220  & 1.47322 & 3.27067 & 772.761 &   314.54       & 488.50 &  1.5531 & 1200.145  \\
    \hline\hline
    \multicolumn{9}{r}{\small{$^\star$ $\Delta V$ not optimized}.}
    \end{tabular}
    }
    \caption{Summary of the best found settlement trees.
    }
    \label{tab:my_label}
\end{table}

\section{Conclusions}
Team Sapienza-PoliTo proposed a solution strategy mainly devoted to maximize the ``uniformity'' of the final settlement tree. In order to reduce the search space dimension, the efforts have been focused on the achievement of an \textit{a priori} chosen star distribution, constructed after a careful analysis of the problem merit function and constraints.
The search procedure has been further simplified by partitioning the overall distribution into smaller, independent, zones, each to be explored departing from one of the root stars encountered by the initial ships. An \textit{ad-hoc} formulated Beam Search procedure has been exploited to built suitable settlement trees in each of the identified zones, leading to an overall, 1013-star, solution, with a distribution quite close to the desired one.
As final step, a concurrent optimization of transfers (by means of an indirect procedure) and transfer times was performed, leading to the final, \nth{7} ranked, submitted solution.
A number of refinement procedures have been developed in the very last days of the competition, but completed immediately after. Thanks to such procedures, a refined 1220-star solution was then obtained, that would (probably) have ranked fifth if it had been submitted until the competition term.

\section{Acknowledgements}
The authors are willing to acknowledge the restless work and the commitment of the team members
here listed in alphabetical order:
Andrea Ferrero, Dario Pastrone, Filippo Masseni and Luigi Mascolo (PoliTo),  Francesco Simeoni (independent) and Danilo Zona (Sapienza).

\bibliographystyle{AAS_publication}   
\bibliography{references}   

\newpage

\end{document}